\documentclass{amsart}

\usepackage{amsmath}
\usepackage{amsfonts}
\usepackage{amssymb}

\numberwithin{equation}{section}

\newtheorem{theorem}[subsection]{Theorem}

\newtheorem{lemma}[subsection]{Lemma}

\newtheorem{conjecture}[subsection]{Conjecture}

\newcommand{\sym}{\text{sym}^2}
\newcommand{\half}{\frac{1}{2}}
\newcommand{\hhalf}{\tfrac{1}{2}}

\newcommand{\h}{\hbar}
\newcommand{\ep}{\epsilon}
\newcommand{\al}{\alpha}
\newcommand{\summ}{\mathop{{\sum}^{h}}}

\begin{document}

\title[Non-vanishing of the symmetric square $L$-function]{Non-vanishing of the symmetric square $L$-function at the central point}
\author{Rizwanur Khan}
\address{University of California at Los Angeles, Department of Mathematics,
         Los Angeles, CA 90095, USA}
\email{rrkhan@ucla.math.edu}

\begin{abstract}
Using the mollifier method, we show that for a positive proportion
of holomorphic Hecke eigenforms of level one and weight bounded by
a large enough constant, the associated symmetric square
$L$-function does not vanish at the central point of its critical
strip. We note that our proportion is the same as that found by
other authors for other families of $L$-functions also having
symplectic symmetry type.
\end{abstract}

\maketitle

\section{Introduction}

A recurring theme in the theory of $L$-functions is the study of
whether or not an $L$-function vanishes at its central point (the
point of symmetry of the $L$-function's functional equation). One
may be interested in this question if the central value is
expected to be a special value, as in the case of the $L$-function
associated to an elliptic curve over a number field, for which the
Birch and Swinnerton-Dyer Conjecture predicts that the order of
vanishing at the central point is the same as the rank of the
elliptic curve. Even if the central value of the $L$-function
under study is not expected to be a special value, the question of
non-vanishing is connected to deep conjectures on the distribution
of the $L$-function's zeros. The following type of problem is
usually formulated: given a suitable collection of related
$L$-functions ordered by analytic conductor, how much of this
family is non-vanishing at the central point? A great deal of work
has been done on this problem for many families of $L$-functions,
of which only a small selection will be referenced in the course
of this paper. Here we study the problem for a family of degree
three $L$-functions: symmetric square $L$-functions lifted from
classical cusp forms.

Let $H_k$ denote the set of holomorphic cusp forms $f$ of level
one and weight $k$ which are eigenfunctions of every Hecke
operator and have Fourier expansions of the form
\begin{align}
f(z)=\sum_{n=1}^{\infty} a_f(n)n^{\frac{k-1}{2}} e^{2\pi i nz}
\end{align}
for $\Im(z)>0$ with $a_f(1)=1$. The coefficients $a_f(n)$ are
necessarily real and satisfy the Ramanujan-Petersson bound
$|a_f(n)|\le d(n)$, where $d(n)$ is the divisor function. The set
$H_k$ has $\frac{k}{12}+O(1)$ elements and forms an orthogonal
basis of the space of cusp forms of level one and weight $k$ (see
\cite{iwa}). For $\Re(s)>1$ we can define the symmetric square
$L$-function associated to $f\in H_k$ as the absolutely convergent
Dirichlet series
\begin{align}
\label{symdiri} L(s,\sym f)= \zeta(2s) \sum_{n=1}^{\infty}
\frac{a_f(n^2)}{n^s}.
\end{align}
In terms of the Rankin-Selberg $L$-function of $f$, we have
$L(s,\sym f)=\frac{\zeta(2s)L(s,f \otimes f)}{\zeta(s)}$. The
symmetric square $L$-function has the degree three Euler product
\begin{align}
\label{euler}L(s,\sym
f)=\prod_{p}(1-a_f(p^2)p^{-s}+a_f(p^2)p^{-2s}-p^{-3s})^{-1}.
\end{align}
Shimura \cite{shim} showed that $L(s,\sym f)$ analytically extends
to an entire function on ${\mathbb C}$ and satisfies the
functional equation
\begin{align}
L_{\infty}(s)L(s,\sym f)=L_{\infty}(1-s)L(1-s,\sym f),
\label{symfunct}
\end{align}
where \begin{align} L_{\infty}(s)=\pi^{-\frac{3}{2}s}
\Gamma(\tfrac{s+1}{2}) \Gamma(\tfrac{s+k-1}{2})
\Gamma(\tfrac{s+k}{2})=\pi^{-\frac{3}{2}s+\half}2^{2-s-k}\Gamma(\tfrac{s+1}{2})\Gamma(s+k-1).
\end{align}
Thus under our normalization the central point is $s=\half$. We
also see from the functional equation that the conductor of
$L(s,\sym f)$ in the weight aspect is $k^2$. Gelbart and Jacquet
\cite{gelbjac} showed that $L(s,\sym f)$ is actually the
$L$-function attached to a cusp form on $GL(3)$.

There is no trivial reason for $L(s,\sym f)$ to vanish at
$s=\half$; that is, the functional equation does not force
$L(\half,\sym f)$ to equal zero. There seems to be no non-trivial
reason either, as $L(\half,\sym f)$ is not known to be a special
value. By this general principle and by the Density Conjecture,
which we discuss in the next section, it is conjectured that
$L(\half,\sym f)$ is never zero. Indeed it would be a great
achievement to prove that for a positive percentage of the family
$H_k$ we have $L(\half,\sym f)\neq 0$. We do not solve this
problem here but settle for something weaker. We show a positive
percentage of non-vanishing on a larger family: essentially
$\bigcup_{\substack{K \le k \le 2K\\k\equiv 0 \bmod 2}} H_k$ for
large $K$. We are in the standard situation where the ratio of
$\log$(conductor) to $\log$(family size) is 1, while this ratio is
2 in the case of a family of an individual weight. In a different
direction, Blomer \cite{blomer} has considered symmetric square
lifts of cusp forms $f$ of prime level $N$, real primitive
nebentypus and fixed weight (so that the ratio
$\log$(conductor):$\log$(family) is $1$). He showed that for large
$N$, we have $L(\half,\sym f)\neq 0$ for a positive proportion of
this family. However he did not compute an explicit proportion.

For any complex numbers $\alpha_f$ depending on $f$, let
\begin{align}
\summ_{f\in H_k} \alpha_f = \frac{2\pi^2}{k-1}\sum_{f\in H_k}
\frac{\alpha_f}{L(1,\sym f)}.
\end{align}
denote the `harmonic' average of $\alpha_f$. The weights $L(1,\sym
f)^{-1}$ arise naturally from the Petersson trace formula. Note
that $\summ_{f\in H_k} 1= 1 + O(2^{-k})$ by the trace formula. Our
main theorem is

\begin{theorem} \label{main}
Let $h$ be a positive valued, infinitely differentiable and
compactly supported function on the positive reals. We have for
any $0<a< \half$ fixed and large enough $K$,
\begin{align*}
\sum_{k \equiv 0 \bmod 2}
h\Big(\frac{k}{K}\Big)\summ_{\substack{f\in H_k\\L(\half,\sym
f)\neq 0}}1 \ge \Big(1-\frac{1}{(1+a)^{3}}\Big)\sum_{k \equiv 0
\bmod 2} h\Big(\frac{k}{K}\Big)\summ_{f\in H_k} 1.
\end{align*}
\end{theorem}
\noindent For the rest of the paper, let us fix a function $h$
satisfying the conditions of the theorem. Thus taking $a$
arbitrarily close to $\half$ and counting harmonically, we have
that at least $\frac{19}{27}>70 \%$ of our family is non-vanishing
at the central point. The reason for the theorem's presentation of
the constant of proportionality in terms of a parameter $a$ rather
than just the best numerical value will be explained.

We expect the same theorem to hold if we replace the sums $\summ$ by
uniform sums $\frac{1}{|H_k|}\sum$. We have only worked with
harmonic sums here, but it would be interesting to remove the
weights $L(1,\sym f)^{-1}$. A method to do this was described in
\cite{kowmic2}. With less work, one can get a smaller but explicit
uniform proportion of non-vanishing from our result using just the
Cauchy-Schwarz inequality and moments of $L(1,\sym f)$.

\subsection{The Katz-Sarnak Philosophy}

In this section, we discuss a few instances relevant to our paper
of the so-called Katz-Sarnak philosophy. Since the discovery of a
connection between zeros of the Riemann Zeta function and
eigenvalues of random matrices by Montgomery \cite{mont} and Dyson
\cite{dyson} in the 1970's, there has been a lot of work by many
authors studying the distribution of zeros of $L$-functions in
families. The conjectures stemming from this work (and the proofs
in function field settings by Katz and Sarnak \cite{katzsarn})
suggest that the distribution of zeros near the central point of a
family of $L$-functions is the same as the distribution of
eigenvalues with small argument of a corresponding classical
compact group. Here we are concerned with the group $USp(2N)$ of
$2N\times 2N$ complex matrices which are both unitary and
symplectic.

In \cite{ils}, Iwaniec, Luo, and Sarnak studied the one-level
density distribution of the low-lying zeros of some $GL(2)$ and
$GL(3)$ families of $L$-functions. Their work yields non-vanishing
results conditional on the Generalized Riemann Hypothesis for
these families. Therefore, let us assume the GRH throughout this
section even though some other results of their paper are
unconditional. For the the symmetric square family $\{ L(s,\sym
f)| f\in H_k \}$, denote a zero of $L(s,\sym f)$ by $\half+i
\gamma_f$, where $\gamma_f$ is real. Let $T$ be a slowly growing
function of $k$, say $T=O(\log k)$. The number of zeros of
$L(s,\sym f)$ counted with multiplicity with imaginary part
between zero and the short height $T$ is asymptotic to $c T \log
k$ for some constant $c$. Let $\tilde{\gamma}_f= \gamma_f c \log k
$ so that the average spacing between these normalized zeros is 1.
Define for an even Schwartz class function $\phi$,
\begin{align}
D(\sym f, \phi)=\sum_{\gamma_f} \phi(\tilde{\gamma}_f).
\end{align}
To get an understanding of the fine behavior of the zeros near the
point $s=\hhalf$, we would like to evaluate the average of $D(\sym
f, \phi)$ over the family $H_k$ for $\phi$ with arbitrarily small
support. Investigating this problem, Iwaniec, Luo, and Sarnak
proved that for $W_{USp}(x)=1-\frac{\sin(2\pi x)}{2\pi x}$, we
have
\begin{align}
\label{ilsdens} \frac{1}{|H_k|}\sum_{f\in H_k} D(\sym f, \phi) \sim
\int_{\mathbb{R}} \phi(x) W_{USp}(x) dx,
\end{align}
as $k\rightarrow \infty$ provided that $\hat{\phi}$, the Fourier
transform of $\phi$, is supported on the interval
$(-\hhalf,\hhalf)$. The condition on $\hat{\phi}$ is very
restrictive and does not permit taking $\phi$ with finite support.
The Density Conjecture claims that (\ref{ilsdens}) holds without
restriction on $\hat{\phi}$. With an extra averaging over the
weight the same result is found for a larger class of test
functions. Assuming the GRH for Dirichlet $L$-functions as well,
the three authors proved that
\begin{align}
\label{ilsdens2} \frac{2}{K} \sum_{\substack{4\le k\le K \\ k \equiv
0 \bmod 2}} \frac{1}{|H_k|} \sum_{f\in H_k} D(\sym f, \phi) \sim
\int_{\mathbb{R}} \phi(x) W_{USp}(x) dx,
\end{align}
as $K\rightarrow \infty$ provided that $\hat{\phi}$ is supported
on the interval $(-\tfrac{3}{2},\tfrac{3}{2})$. Notice that one
expects a repulsion of zeros from the central point: when $|x|\le
\ep$ for a small positive constant $\ep$, we have $1-\frac{\sin
(2\pi x)}{2\pi x}\ll \ep^2$ and so $\int_{-\ep}^{\ep} W_{USp}(x)
dx \ll \ep^3$.

What is amazing is that $W_{USp}$ is the same function which
occurs on the Random Matrix Theory side. A matrix $A\in U(2N)$ has
$N$ pairs of complex conjugate eigenvalues $e^{i \theta_A}$ on the
unit circle, counted with multiplicity. Let $\tilde{\theta}_A=
\frac{N}{\pi} \theta_A$ so that the average spacing between the
normalized angles is 1. Katz and Sarnak showed that for an even
Schwartz class function $\phi$ we have
\begin{align}
\int_{USp(N)} \sum_{\theta_A} \phi\big( \tilde{\theta}_A \big) dA
\sim \int_{\mathbb{R}} \phi(x) W_{USp}(x) dx,
\end{align}
as $N \rightarrow \infty$, where $dA$ is Haar measure.

Another family that is found to be symplectic is the family of
degree one quadratic Dirichlet $L$-functions $
L(s,(\frac{8d}{\cdot}))$ for $d$ odd and squarefree with $X\le d
\le 2X$. This was first shown by \"{O}zluk and Snyder
\cite{snyder}. For another example, consider the set of degree two
$L$-functions attached to `odd' primitive cusp forms $f$ of weight
2 and prime level $q$. By odd we mean that at the central point,
the functional equation of $L(s,f)$ reads
$L(\half,f)=-L(\half,f)$. Thus in the family of such Hecke
$L$-functions, it is trivial that $L(\half,f)=0$. This family is
found to be of type $SO(2N+1)$, with corresponding density
function
\begin{align}
\label{orthsymp} W_{SO(2N+1)}=\delta_0+W_{USp(2N)}.
\end{align}
Above $\delta_0$ is the Dirac delta function, which occurs because
every (unitary) matrix in $SO(2N+1)$ has $1$ as an eigenvalue. In
this case it is clearly more interesting to study the
non-vanishing of the derivative $L^{\prime}(s,f)$.  Since
$L'(s,f)$ vanishes precisely when $L(s,f)$ has an additional zero
at $s=\half$, we have by the relation (\ref{orthsymp}) that the
family of derivatives is of symplectic type.

Let us now apply the density formula (\ref{ilsdens2}) to the
non-vanishing of the symmetric square $L$-function at the central
point. For the simple choice of $\phi(x)=\Big(\frac{\sin (v \pi x
)}{v \pi x}\Big)^2$ where $v>0$, we have $\phi(0)=1$, $\phi(x)\ge
0$, and $\hat{\phi}$ supported on $(-v,v)$. With this test
function we have from (\ref{ilsdens2}) that
\begin{align}
\frac{2}{K} \sum_{\substack{4\le k\le K \\ k \equiv 0 \bmod 2}}
\frac{1}{|H_k|}\sum_{f\in H_k} ord_{s=\half} L(s,\sym f) \le
\int_{\mathbb{R}} \phi(x) W_{USp}(x) dx.
\end{align}
Using the fact that the order of vanishing of $L(s,\sym f)$ at
$s=\half$ must be even from the sign of the functional equation,
it follows that for a proportion of at least
\begin{align}
1-\frac{1}{4v^2} \label{formr}
\end{align}
of the family $\{ L(s,\sym f)|f\in \bigcup_{\substack{K \le k \le
2K\\k\equiv 0 \bmod 2}} H_k \}$ for large enough $K$ we have that
$L(\hhalf, \sym f)$ is nonzero. Thus taking $v=\frac{3}{2}$ the
proportion $\frac{8}{9}$ is gotten conditionally on GRH in
\cite{ils}. Actually a slightly better answer is also found by
optimizing the choice of $\phi$.

Similarly for the quadratic Dirichlet $L$-function family above,
by taking $v=2$ in (\ref{formr}), \"{O}zluk and Snyder
conditionally showed that for large $X$ at least $\frac{15}{16}$
of the family is non-vanishing at the central point. For the third
family, of odd Hecke $L$-functions, again it was conditionally
obtained in \cite{ils} that when $q$ is large enough, for at least
$\frac{15}{16}$ of the family we have $L^{\prime}(f,\hhalf)\neq
0$. In all three examples, the Density Conjecture implies a
proportion of 1 by taking $v$ arbitrarily large.

\subsection{The mollifier method}

We prove Theorem \ref{main} using the mollifier method, a
technique which goes back to Bohr and Landau and was used by
Selberg \cite{selb} to show that a positive proportion of the
non-trivial zeros of the Riemann Zeta function lie on the half
line. Let us begin by observing that in principle we can recover
the distribution of the values $L(\half,\sym f)$ for $f\in H_k$
from asymptotics for $\summ_{f\in H_k} L(\half,\sym f)^n$ for all
$n\in {\mathbb N}$. However our knowledge of these moments is very
limited. For the first moment it is known (see \cite{lau}) that
\begin{align}
\summ_{f\in H_k} L(\hhalf,\sym f)\sim \log k.
\end{align}
For the second
moment we have the conjecture
\begin{conjecture}\label{2conj}
For some constant $c$,
\begin{align}
\summ_{f\in H_k}L(\hhalf,\sym f)^2 \sim c (\log k)^3.
\end{align}
\end{conjecture}
\noindent Even the upper bound $\summ_{f\in H_k}L(\hhalf,\sym
f)^2\ll k^{\epsilon}$ is not established. The difficulty of this
conjecture is what brings us to take an extra averaging over the
weight. We shall show
\begin{align}
&\nonumber\sum_{k \equiv 0 \bmod 2} h\Big(\frac{k}{K}\Big)
\summ_{f\in H_k} L(\hhalf,\sym f) = K P_h(\log K)+ O(K^{\ep}),\\
&\label{simplemoments} \sum_{k \equiv 0 \bmod 2}
h\Big(\frac{k}{K}\Big) \summ_{f\in H_k} L(\hhalf,\sym f)^2 = K
Q_h(\log K)+O(K^{\ep}),
\end{align}
where $P_h$ and $Q_h$ are degree one and degree three polynomials
respectively which depend on $h$, and $\epsilon$ is an arbitrarily
small positive constant. This shows that the average value of
$L(\hhalf,\sym f)$ is proportional to $\log K$, and we would like
to be able to conclude that a lot of values of $L(\hhalf,\sym f)$
are not zero. However in this way we cannot rule out the
possibility of many values of $L(\half,\sym f)$ being zero and
some being very large. In fact, a comparison of the main terms of
the first and second moments shows that there are fluctuations in
the size of $L(\hhalf,\sym f)$. Nevertheless by the Cauchy-Schwarz
inequality we get that
\begin{align}
\sum_{k \equiv 0 \bmod 2}
h\Big(\frac{k}{K}\Big)\summ_{\substack{f\in H_k\\L(\half,\sym
f)\neq 0}}1&\ge \frac{(\sum_{k \equiv 0 \bmod 2}
h(\tfrac{k}{K})\summ_{f\in H_k} L(\hhalf,\sym f))^2}{\sum_{k
\equiv 0 \bmod 2} h(\tfrac{k}{K})\summ_{f\in
H_k}L(\hhalf,\sym f)^2}\\
&\nonumber \gg (\log K)^{-1}K.
\end{align}
Thus we get that for a proportion of $(\log K)^{-1}$ of Hecke cusp
forms of weight less than $K$, we have $L(\hhalf,\sym f)\neq 0$.
This however is $0 \%$ and we can improve upon it by using the
mollifier method. With the power savings in the error terms of
(\ref{simplemoments}), there is room to find a little more than
the second moment. Define a short Dirichlet series, called a
mollifier,
\begin{align}
M(f) = \sum_{r\le K^a} \frac{x_r a_f(r^2)}{\sqrt{r}},
\end{align}
where $a>0$ is a constant and $x_r$ are coefficients to be chosen.
Since the non-vanishing of $L(\half,\sym f)M(f)$ implies the
non-vanishing of $L(\half,\sym f)$, we have as before,
\begin{align}
\label{ratiointro} \sum_{k \equiv 0 \bmod 2}
h\Big(\frac{k}{K}\Big)\summ_{\substack{f\in H_k\\L(\half,\sym
f)\neq 0}}1 \ge \frac{(\sum_{k \equiv 0 \bmod 2}
h(\tfrac{k}{K})\summ_{f\in H_k} L(\hhalf,\sym f)M(f))^2}{\sum_{k
\equiv 0 \bmod 2} h(\tfrac{k}{K})\summ_{f\in H_k} L(\hhalf,\sym
f)^2M(f)^2}.
\end{align}
We shall compute the mollified moments in terms of the
coefficients $x_r$ and then carefully choose $x_r$ to maximize
this ratio. This choice of $M(f)$ dampens or mollifies the large
values of $L(\hhalf,\sym f)$, so that the square of the first
moment and $K$ times the second moment are comparable. Thus we
will get that $L(\hhalf,\sym f)$ is not zero for a positive
percentage of our family. In fact we will show that for an optimal
mollifier of length $K^{a}$, where recall that $K^2$ is the size
of conductor of our family, we get a proportion of $1-(1+a)^{-3}$.
Our methods enable us to take $a<\half$.

Soundararajan \cite{soun} proved using the mollifier method that
for at least $\frac{7}{8}$ of positive odd square-free integers
$X\le d\le 2X$ we have $L(\half,(\frac{8d}{\cdot}))\neq 0$ when
$X$ is large enough. Taking a mollifier of length $(\sqrt{X})^a$,
where $X$ is the size of conductor of any quadratic Dirichlet
$L$-function in this family, he obtained a proportion of
$1-(1+a)^{-3}$ and was able to take $a<1$ to get his result.
Kowalski and Michel \cite{kowmic} studied the non-vanishing of
$L^{\prime}(f,\hhalf)$ for odd primitive cusp forms $f$ of weight
2 and prime level $q$. Taking a mollifier of length
$(\sqrt{q})^a$, where $q$ is the conductor of $L'(s,f)$ and
$a<\frac{1}{2}$ (which the authors remarked can be improved to
$a<1$), they showed that for large enough $q$ we have that
$L^{\prime}(f,\hhalf)\neq 0$ for a proportion of $1-(1+a)^{-3}$ of
this family. Thus our result supports the belief that these three
different families of $L$-functions share the same distribution of
low lying zeros. In fact, from a family's symmetry type one can
make predictions using the Ratios Conjecture about non-vanishing
by the use of mollifiers (see \cite{conrsnai}).

Being able to take a mollifier $M(f)$ of length $K^a$ is roughly
comparable to being able find the $(2+2a)^{\text{th}}$ moment of
$L(\half,\sym f)$ on average over the weight. Taking $\hat{\phi}$
of support $(-v,v)$ in the density formula (\ref{ilsdens2}) is
comparable to finding the $(2v)^{\text{th}}$ moment of
$L(\half,\sym f)$ on average over the weight. Thus we can make a
heuristic connection between $a$ and $v$, with $v$ corresponding
to $a+1$. In this way, our proportion of non-vanishing of
$\tfrac{19}{27}$ is the unconditional analogue of the conditional
proportion $\tfrac{8}{9}$ of \cite{ils}. Similarly, the
unconditional proportion $\tfrac{7}{8}$ of the other symplectic
families considered above compares with the conditional proportion
$\tfrac{15}{16}$. Notice also that the proportions $\frac{7}{8}$
and $\frac{19}{27}$ are quite good- this can be explained by the
repulsion of zeros that is expected at the central point in the
symplectic family. This may be compared with the orthogonal family
considered in \cite{sarniwan2}, where a proportion greater than
$\half$ is sought, but not quite achieved, to show that
Landau-Siegel zeros do not exist for Dirichlet $L$-functions.

\section{The first twisted moment}

In this section we find the first moment of $L(\half,\sym f)$
twisted by a Fourier coeffiecient $a_f(r^2)$ on average over $k$.
This will yield the first mollified moment.

\begin{theorem}\label{M1} For $r\le K^{2-\ep}$ we have
\begin{align*}
\sum_{k \equiv 0 \bmod 2} h\Big(\frac{k-1}{K}\Big) \summ_{f \in
H_k}
 L(\hhalf,\sym f)a_f(r^2)=&\frac{K}{2\sqrt{r}}\int_{0}^{\infty} h(u)
\Big(\log(uK/r)+C\Big)du\\
&+O(r^{\half}K^{\ep}),
\end{align*}
where $C$ is an absolute constant, $\epsilon>0$ is an arbitrarily
small positive constant and the error term depends on $\epsilon$ and
$h$.
\end{theorem}
\noindent Let us adopt the following convention throughout this
paper: $\ep$ will always denote an arbitrarily small positive
constant, but not necessarily the same one from one occurrence to
the next. Any implied constants may depend implicitly on $h$ and
$\ep$, unless otherwise indicated. Before proving this theorem we
will need some preliminary results.

\subsection{Petersson trace formula}

We will need the Petersson trace formula:
\begin{align}
\label{petersson} \summ_{f\in H_k} a_{f}(n) a_f(m) =
\delta_{m,n}+2\pi i^k \sum_{c=1}^{\infty}
\frac{S(n,m;c)}{c}J_{k-1}\Big(\frac{4\pi\sqrt{mn}}{c}\Big),
\end{align}
where the value of $\delta_{m,n}$ is $1$ if $m=n$ and $0$ otherwise,
$S(n,m;c)$ is a Kloosterman sum, and $J_{k-1}$ is the $J$-Bessel
function.

\subsection{Approximate functional equation}

While the Dirichlet series expansion of $L(s,\sym f)$ given in
(\ref{symdiri}) is only valid for $\Re(s)>1$, we can use a
standard tool of analytic number theory called an approximate
functional equation to express $L(\half,\sym f)$ as a weighted
Dirichlet series. From property (\ref{V1}) below, the terms of
this series are only significant for $n\le k^{1+\ep}$, so that the
length of the Dirichlet series is essentially equal to the square
root of the conductor. We will use Stirling's approximation: $\log
\Gamma(z)=(z-\hhalf)\log
z-z+\half\log\sqrt{2\pi}+\sum_{n=1}^{m}c_n z^{-2n+1}+O_m(z^{-m})$
when $|\arg z|<\pi-\ep$ for some constants $c_n$.

\begin{lemma} {\bf Approximate functional equation} \label{approxeq}
We have
\begin{align}
L(\hhalf,\sym f) = 2\sum_{n\ge 1} \frac{a_f(n^2)}{n^{\half}}
V_{k}(n), \label{mainapproxfuncteq}
\end{align}
where for any real $\xi>0$ and $A>0$,
\begin{align}
V_{k}(\xi)=\frac{1}{2\pi i} \int_{(A)}
\frac{L_{\infty}(\half+y)}{L_{\infty}(\half)} \zeta(1+2y) \xi^{-y}
\frac{dy}{y}
\end{align}
is real valued and satisfies,
\begin{align}
&V_{k}(\xi) \ll_A \Big(\frac{k}{\xi}\Big)^A, \text{ for any }A>0. \label{V1} \\
&V_{k}^{(B)} (\xi) \ll_{A,B} \xi^{-B} \Big(\frac{k}{\xi}\Big)^A, \text{ for any }A>0 \text{ and integer }B\ge 0. \label{V2} \\
&V_{k}(\xi)= \hhalf(\log(k/\xi)+C)+O\Big(\frac{\xi}{k}\Big), \text{
where } C=2\gamma-\frac{3\log \pi}{2}-\log
2+\frac{\Gamma^{\prime}(\frac{3}{4})}{2\Gamma(\frac{3}{4})}.
\label{V4}
\end{align}
We also have that
\begin{align} V_k(\xi) =\frac{1}{2\pi i} \int_{(A)}
G(y)
\frac{\zeta(1+2y)}{y}\Big(\frac{k}{\xi}\Big)^{y}dy+O(\xi^{-\ep}k^{-1+\ep}),\label{V5}
\end{align}
for any $A>0$, where $G(y)=
\pi^{-\frac{3}{2}y}2^{-y}\frac{\Gamma(\frac{y}{2}+\frac{3}{4})}{\Gamma(\frac{3}{4})}$
exponentially decays on vertical lines.
\end{lemma}

\proof

Define for $-\frac{1}{2} \le c \le \frac{3}{2}$,
\begin{align}
 I(c)=\frac{1}{2\pi i} \int_{(c)} L(\hhalf+y,\sym f)\frac{L_{\infty}(\hhalf+y)}{L_{\infty}(\half)}\frac{dy}{y},
\end{align}
where the integral is taken from $c-i\infty$ to $c+i\infty$ and
converges absolutely because the integrand decays exponentially as
$|\Im(y)|\rightarrow \infty$. In the range $\Re(y)\ge -\half$ the
integrand has a simple pole at $y=0$ with residue $L(\half,\sym
f)$. Thus by Cauchy's theorem we have
$I(\tfrac{3}{2})-I(-\hhalf)=L(\half,\sym f)$. Using the functional
equation (\ref{symfunct}) and a change of variables we have
$I(-\hhalf)=- I(\hhalf)$ and so we get that $L(\half,\sym f)=
I(\tfrac{3}{2})+I(\hhalf)=2I(\tfrac{3}{2})$. At
$\Re(y)=\frac{3}{2}$ we can use the Dirichlet series expansion
(\ref{symdiri}) to get
\begin{align}
 I(\tfrac{3}{2})=\sum_{n\ge 1}
\frac{a(n^2)}{n^{\half}} \int_{(3/2)}
\frac{\zeta(2(\hhalf+y))L_{\infty}(\half+y)}{n^y
L_{\infty}(\half)}\frac{dy}{y}=\sum_{n\ge 1}
\frac{a(n^2)}{n^{\half}} V_{k}(n),
\end{align}
where we implicity exchanged summation and integration by absolute
convergence. This establishes (\ref{mainapproxfuncteq}).

For $\Re(y)=A$ we have by Stirling's estimates that
\begin{align}
\Big|\frac{L_{\infty}(\half+y)}{L_{\infty}(\half)}\Big|\ll
|\Gamma(\tfrac{3+2y}{4})|\frac{\Gamma(k+A-\half)}{\Gamma(k-\half)}
\ll_A |\Gamma(\tfrac{3+2y}{4})|k^A.
\end{align}
\noindent This implies (\ref{V1}). We obtain (\ref{V2}) by first
differentiating $V_k(\xi)$ and then using Stirling's estimates. To
get $(\ref{V4})$ we move the line of integration of $V_k(\xi)$ to
$\Re(y)=-1$, crossing a double pole at $y=0$. The main term
$(\ref{V4})$ is the residue from this pole and the error is given by
Stirling's estimates.

Let us turn to (\ref{V5}). We have
$\tfrac{L_{\infty}(\half+y)}{L_{\infty}(\half)}=G(y)\frac{\Gamma(y+k-\half)}{\Gamma(k-\half)}$.
For $\Im(y)\le k^\ep$ it follows from Stirling's estimates that
$\frac{\Gamma(y+k-\half)}{\Gamma(k-\half)}=k^y(1+O_A(k^{-1+\ep}))=k^y+O_A(k^{A-1+\ep}).$
Since
$\Big|\frac{\Gamma(y+k-\half)}{\Gamma(k-\half)}\Big|\le\frac{\Gamma(A+k-\half)}{\Gamma(k-\half)}\ll_A
k^A$ and $G(y)$ is exponentially decreasing on vertical strips, we
can restrict the integral in $V_k(\xi)$ to $\Im(y)\le k^{\ep}$
with an error of $O_A((k/\xi)^{A}e^{-k^{\ep}})$. Taking $A=\ep$ we
get $V_k(\xi) =\frac{1}{2\pi i} \int_{(\ep)} G(y)
\frac{\zeta(1+2y)}{y}\Big(\frac{k}{\xi}\Big)^{y}dy+O(\xi^{-\ep}
k^{-1+\ep})$. The line of integration can be moved any $A>0$ to
get (\ref{V5}).
\endproof

\subsection{An average of the $J$-Bessel function}

We have the following standard estimates: $|J_{k-1}(t)|\le 1$ for
all $t\ge 0$, $J_{k-1}(t)\sim \frac{1}{\Gamma(k)}(t/2)^{k-1}$ for
$t\le k^{1/2-\ep}$ and $J_ {k-1}(t)\sim (\frac{\pi}{2} t)^{-1/2}
\cos (t-\frac{\pi}{2} k +\frac{\pi}{4})$ for $t\ge k^{2+\ep}$. The
long term behavior of the $J$-Bessel function manifests itself
when we average over $k$:

\begin{lemma}\label{avg}
We have for $t> 0$,
\begin{align}
\label{avgmaint} \sum_{k \equiv 0 \bmod 2} 2i^k
h\Big(\frac{k-1}{K}\Big)
&J_{k-1}(t)\\
&\nonumber = -\frac{K}{\sqrt{t}} \Im \Big(e^{-2\pi i /8} e^{it}
\h\Big(\frac{K^2}{2t}\Big) \Big) +
O\Big(\frac{t}{K^4}\int_{-\infty}^{\infty} v^4|\hat{h}(v)|dv\Big),
\end{align}
where $\h(v)=\int_{0}^{\infty} \frac{h(\sqrt{u})}{\sqrt{2\pi u}}
e^{iuv} du$ and $\hat{h}$ denotes the Fourier transform of $h$.
The implied constant is absolute.
\end{lemma}
\proof We refer to \cite[Lemma 5.8]{iwa} for the full proof. It is
shown there that the sum we want equals
\begin{align}
&2\sum_{k \equiv 0 \bmod 4} h\big(\tfrac{k-1}{K}\big)
J_{k-1}(t)-2\sum_{k \equiv 2 \bmod 4} h\big(\tfrac{k-1}{K}\big)
J_{k-1}(t)\\
&\nonumber=-2\int_{-\infty}^{\infty} K \hat{h}(Kv)  \sin (t\cos
2\pi v)
dv\\
&\nonumber=-2\int_{-\infty}^{\infty} K \hat{h}(Kv)  \sin
(t-2\pi^2t v^2)dv +
O\Big(\int_{-\infty}^{\infty}tv^4K|\hat{h}(Kv)|dv\Big),
\end{align}
where the last line follows by writing $\cos 2\pi v = 1 - 2\pi^2
v^2+O(v^4)$. We have written the main term as it appears in
\cite[Corollary 8.2]{ils}.
\endproof

\noindent We have that $\h(v)\ll 1$ and by integrating by parts
several times we get that $\h(v)\ll_B v^{-B}$
 for any $B\ge0$. Thus the main
term of (\ref{avgmaint}) is not dominant if $t \le K^{2-\ep}$.

Define
\begin{align}
W_K(n,m,v)=\int_{0}^{\infty}
\frac{V_{\sqrt{u}K+1}(n)V_{\sqrt{u}K+1}(m)h(\sqrt{u})}{\sqrt{2\pi
u}}e^{iuv} du.
\end{align}
By lemma \ref{avg} we have
\begin{align}
\label{avgwithV} \sum_{k \equiv 0 \bmod 2} 2i^k
h\Big(\frac{k-1}{K}\Big) &V_k(n)V_k(m)
J_{k-1}(t)=-\frac{K}{\sqrt{t}}
\Im \Big(e^{-2\pi i /8} e^{it} W_K(n,m,\tfrac{K^2}{2t}) \Big) \\
&\nonumber+ O\Big(\frac{t}{K^4} \int_{-\infty}^{\infty} v^4\Big|
\int_{0}^{\infty} V_{uK+1}(n)V_{uK+1}(m)h(u)e^{iuv} du \Big|
dv\Big).
\end{align}
Using the integral definition of $V_k(n)$ we have
\begin{align}
 W_K(n,m,v)= &\frac{1}{(2\pi i)^2} \int_{(A_1)}\int_{(A_2)}
G(y)G(x) \frac{\zeta(1+2y)}{y}\frac{\zeta(1+2x)}{x}\frac{1}{m^xn^y}
\\
&\int_{0}^{\infty}
\frac{\Gamma(\sqrt{u}K+y+\half)}{\Gamma(\sqrt{u}K+\half)}\frac{\Gamma(\sqrt{u}K+x+\half)}{\Gamma(\sqrt{u}K+\half)}
\frac{h(\sqrt{u})}{\sqrt{2\pi u}} e^{iuv} du dy dx\nonumber,
\end{align}
for any $A_1,A_2>0$. By integrating by parts several times and
Stirling's formula we have for $\Re(y)=A_2$, $|y|\le K^{\ep}$,
$\Re(x)=A_1$ and $|x|\le K^{\ep}$ that
\begin{align}
\int_{0}^{\infty}
\frac{\Gamma(\sqrt{u}K+y+\half)}{\Gamma(\sqrt{u}K+\half)}\frac{\Gamma(\sqrt{u}K+x+\half)}{\Gamma(\sqrt{u}K+\half)}
&\frac{h(\sqrt{u})}{\sqrt{2\pi u}} e^{iuv} du\\
&\nonumber \ll_{B,A_1,A_2}\frac{(1+|x|+|y|)^B K^{A_1+A_2}}{v^B},
\end{align}
for any $B\ge 0$. This together with the fact that $G(y)$
decreases exponentially as $\Im(y)\rightarrow \infty$ implies
\begin{align}
\label{wbound} W_K(n,m,v)\ll_{B,A_1,A_2} (\tfrac{K}{n})^{A_1}
(\tfrac{K}{m})^{A_2}v^{-B}.
\end{align}
Thus $W$ is essentially supported on $n\le K^{1+\ep}$, $m\le
K^{1+\ep}$, and $v\le K^{\ep}$. Similarly by integrating by parts,
the error in (\ref{avgwithV}) is
\begin{align}
O_{A_1,A_2}\Big(\frac{t}{K^4}
\Big(\frac{K}{n}\Big)^{A_1}\Big(\frac{K}{m}\Big)^{A_2}\Big),
\end{align}
for any $A_1,A_2>0$. In applications $t$ will be bounded by a
power of $K$ and so this error term essentially only appears when
$n\le K^{1+\ep}$ and $m\le K^{1+\ep}$.

In combination with (\ref{V5}) we have
\begin{align}
\label{wtow}
W_K(n,m,v)=W(\tfrac{n}{K},\tfrac{m}{K},v)+O(K^{-1+\ep}),
\end{align}
where we define for real $\xi_1>0$ and $\xi_2>0$,
\begin{align}
\label{definw} W(\xi_1,\xi_2,v)=\frac{1}{(2\pi i)^2}
\int_{(A_1)}\int_{(A_2)} G(y)G(x)
\frac{\zeta(1+2y)}{y}\frac{\zeta(1+2x)}{x} \xi_1^{-y} \xi_2^{-x}
\h_{x+y}(v) dy dx
\end{align}
and
\begin{align}
\label{hbardef} \h_{z}(v)=\int_{0}^{\infty}
\frac{h(\sqrt{u})}{\sqrt{2\pi u}} u^{z/2} e^{iuv} du,
\end{align}
for a complex number $z$. By integration by parts we have
\begin{align}
\label{hbarbound} \h_{z}(v)\ll_{\Re(z),B} (1+|z|)^B v^{-B}
\end{align}
and
\begin{align}
\label{w1bound} W(\xi_1,\xi_2,v)\ll_{B,A_1,A_2}
\xi_1^{-A_1}\xi_2^{-A_2} v^{-B},
\end{align}
for $A_1,A_2>0$ and $B\ge 0$.

Similarly we have
\begin{align}
\label{avg3} \sum_{k \equiv 0 \bmod 2} 2i^k
h\Big(\frac{k-1}{K}\Big) V_k(n) J_{k-1}(t)=&-\frac{K}{\sqrt{t}}
\Im \Big(e^{-2\pi i /8} e^{it}
W_K(n,\tfrac{K^2}{2t}) \Big)\\
&\nonumber \indent +O_A\Big(\frac{t}{K^4}
\Big(\frac{K}{n}\Big)^{A} \Big),
\end{align}
where $W_K(n,v)=\int_{0}^{\infty}
\frac{V_{\sqrt{u}K+1}(n)h(\sqrt{u})}{\sqrt{2\pi u}} e^{iuv} du$ and
$A>0$. We have
\begin{align}
\label{wnbound} W_K(n,v)\ll_{A,B} (\tfrac{K}{n})^A v^{-B},
\end{align}
for any $A>0$ and $B\ge 0$.

\subsection{Kloosterman sums}

\noindent We will use Weil's bound for the Kloosterman sum,
\begin{align}
|S(n,m,c)|\le d(c)c^{\half}(n,m,c)^{\half},
\end{align}
and the average bound which follows from this,
\begin{align}
\sum_{n\le x} |S(n^2,m^2,c)|\ll c^{\half+\ep} \sum_{n\le x}
(n^2,c)^{\half}\ll c^{\half+\ep}
\sum_{\alpha|c}\sum_{\substack{n\le x\\(n,c)=\alpha}}\alpha \ll x
c^{\half+\ep}.
\end{align}

\subsection{Proof of Theorem \ref{M1}}

Using the approximate functional equation and the Petersson trace
formula we get that the first twisted moment is
\begin{align}
&\label{first} \sum_{k \equiv 0 \bmod 2} h\Big(\frac{k-1}{K}\Big)
\summ_{f \in H_k}
 L(\hhalf,\sym f) a_f(r^2) \\
&\indent=2\sum_{k \equiv 0 \bmod 2} h\Big(\frac{k-1}{K}\Big)
\Big(\frac{V_k(r)}{\sqrt{r}}+\sum_{n,c\ge 1} 2 \pi i^k
\frac{S(n^2,r^2,c)}{c}J_{k-1}\Big(\frac{4\pi nr}{c}\Big)
\frac{V_k(n)}{\sqrt{n}}\Big)\nonumber.
\end{align}
Estimating $V_k(r)$ using (\ref{V4}) and using (\ref{avg3}) we get
that (\ref{first}) equals,
\begin{align}
&\sum_{k \equiv 0 \bmod 2} h\Big(\frac{k-1}{K}\Big)\frac{
\log(k/r)+C+O(r/K)}{\sqrt{r}}\label{er3}\\
&\nonumber -2\pi\sum_{n \ge 1}\frac{1}{\sqrt{n}}K\sum_{c\ge
1}\Im\Big(e(-1/8)\frac{1}{\sqrt{4\pi
nr/c}}e(2nr/c)\frac{S(n^2,r^2,c)}{c}W_K(n,\tfrac{K^2c}{8\pi
nr})\Big) \\
&\nonumber\hspace{1.9in}+ O\Big(K^{-4+\ep}\sum_{n,c\ge
1}r\sqrt{n}\frac{|S(n^2,r^2,c)|}{c^2}\Big(\frac{K}{n}\Big)^A\Big),
\end{align}
for any $A>0$. In the third line of (\ref{er3}), the sum can be
restricted to $n\le K^{1+\ep}$ with an error of $O(K^{-10})$. By
the bound (\ref{wnbound}) on $W_K$, the sums in the second line
can be restricted to $n\le K^{1+\ep}$ and $c\le nr/K^{2-\ep}$,
with an error of $O(K^{-10})$. Then using Weil's bound for the
Kloosterman sum, we find that the last two lines of $(\ref{er3})$
are bounded by $O(r^{\half}K^{\ep})$. The error in the first line
of (\ref{er3}) is also $O(r^{\half}K^{\ep})$. Finally for the main
term, note that by Poisson summation we have
\begin{align}
\sum_{k \equiv 0 \bmod 2}
h\Big(\frac{k-1}{K}\Big)=\frac{K}{2}\hat{h}(0)+
O_B(K^{-B})\label{hsum},
\end{align}
for any $B>0$, where $\hat{h}$ denotes the Fourier transform of $h$.
Writing $\log k=\log(\tfrac{k-1}{K}K)+O(k^{-1})$ and using
(\ref{hsum}) the proof is complete.

\section{The  second twisted moment}

In this section we find the second moment of $L(\half,\sym f)$
twisted by a Fourier coeffiecient $a_f(r^2)$ on average over $k$.
This will yield the second mollified moment. Interesting features
make this computation natural and simple.

\begin{theorem}\label{M2}
For any positive integer $d$, write $d=d_1 d_2^2$ with $d_1$
square-free. For $r\le K^{1-\ep}$ we have
\begin{align*}
&\sum_{k \equiv 0 \bmod 2} h\Big(\frac{k-1}{K}\Big) \summ_{f \in
H_k}
L(\hhalf,\sym f)^2a_f(r^2)=\\
\nonumber&\hspace{0.3in} \frac{K}{\sqrt{r}}\sum_{d|r^2}
\frac{1}{\sqrt d_1} \int_{0}^{\infty} h(u)
\Big(\frac{1}{4}\log^2\frac{uK}{d_1d_2}\log\frac{uKd_2}{r}-\frac{1}{12}\log^3\frac{uK}{d_1d_2}
+P_2\Big(\log\frac{uK}{d_1d_2}\Big)\\
&\nonumber \hspace{2.4in}
+\log\frac{uKd_2}{r}P_1\Big(\log\frac{uK}{d_1d_2}\Big)
\Big)du+O(r^{\half}K^{\ep}).
\end{align*}
where $P_i$ is a polynomial of degree $i$ and we write $\log^i x =
(\log x)^i$.
\end{theorem}
\noindent We will need the following results.

\subsection{Sum of Kloosterman sums}

We will need the following sum of twisted Kloosterman sums.

\begin{lemma}\label{sumkloo}
Suppose $n,b,d$ and $r$ are positive integers with $r\le K^{1-\ep}$.
We have
\begin{align}
\label{kloosum}&\sum_{n\ge 1} S(n^2,b^2,c)e(2nb/c)n^{-1}W(\tfrac{n}{K},\tfrac{bd}{rK},\tfrac{K^2c}{8\pi n b})\\
&\nonumber\indent=\begin{cases} \phi(c) c^{-\half} \sum_{n\ge
1}n^{-1}W(\tfrac{n}{K},\tfrac{bd}{rK},\tfrac{K^2c}{8\pi n b })+
O(K^{-100}) &\text{ if }c \text{ is a square,}\\
O(K^{-100}) &\text{ otherwise}.
\end{cases}
\end{align}
\end{lemma}

\proof We have
\begin{align}
&\nonumber \sum_{n\ge 1} S(n^2,{b}^2,c)e(2nb/c)n^{-1}W(\tfrac{n}{K},\tfrac{bd}{rK},\tfrac{K^2c}{8\pi nb})\\
&\label{splitresidues} \hspace{1.3in} =\sum_{a\bmod c}
S(a^2,b^2,c)e(2ab/c) \sum_{\substack{n\ge 1\\ n\equiv a \bmod c}}
n^{-1}W(\tfrac{n}{K},\tfrac{bd}{rK},\tfrac{K^2c}{8\pi nb}).
\end{align}

We first deal with the inner sum above. We assume throughout this
proof that $b\le rK^{1+\ep/4}\le K^{2-\ep/2}$ and $c\le nb
K^{-2+\epsilon/4} \le  K^{1-\ep/2}$, since otherwise both sides of
(\ref{kloosum}) would be of size $O(K^{-100})$. This follows by
the bound (\ref{w1bound}). Using definition (\ref{definw}) and
integration by parts we have have the following bound on the
derivatives of $W$ for $\xi>0$:
\begin{align}
\label{wderivboundw} &\frac{\partial^{j}}{\partial \xi^{j}}
W(\xi,\tfrac{bd}{rK},\tfrac{Kc}{8\pi \xi b}) \ll_{j,A,B}
\xi^{-j-A}(\tfrac{\xi b}{K c})^{B}K^{\ep},
\end{align}
for any $A > 0$, $B\ge0$ and integer $j\ge 0$. We extend the
definition of $W(\xi,\tfrac{bd}{rK},\tfrac{Kc}{8\pi \xi b})$ to
all real numbers $\xi$ by setting
$W(\xi,\tfrac{bd}{rK},\tfrac{Kc}{8\pi \xi b})=0$ for $\xi \le0$;
by the bound just shown this is still a smooth function of $\xi$.
Now by Poisson summation we have that the inner sum of
(\ref{splitresidues}) equals
\begin{align}
\nonumber \sum_{ n\equiv a \bmod c}
&n^{-1}W(\tfrac{n}{K},\tfrac{bd}{rK},\tfrac{K^2c}{8\pi nb})\\
&\label{smoothsum} \hspace{1in}= \frac{1}{c} \sum_{n}
\int_{-\infty}^{\infty} \xi^{-1}
W(\xi,\tfrac{bd}{rK},\tfrac{Kc}{8\pi \xi b}) e^{ \frac{2\pi i
K}{c} \xi n}e^{ \frac{-2\pi i a}{c} n} d\xi.
\end{align}
By integrating by parts $j$ times and using (\ref{wderivboundw})
we find that for $|n|\neq 0$ we have
\begin{align}
\int_{-\infty}^{\infty} \xi^{-1}
W(\xi,\tfrac{bd}{rK},\tfrac{Kc}{8\pi \xi b}) e^{ \frac{2\pi i
K}{c} \xi n}e^{ \frac{-2\pi i a}{c} n} d\xi \ll_j
(\tfrac{c}{K|n|})^{j}(\tfrac{ b}{K c})^{j}K^{10}.
\end{align}
Since $b\le K^{2-\ep/2}$, the above is of size
$O(|n|^{-2}K^{-100})$ by taking $j$ large enough. Thus only the
$n=0$ term is significant and (\ref{smoothsum}) equals
\begin{align}
 \frac{1}{c} \int_{-\infty}^{\infty}  \xi^{-1}
W(\xi,\tfrac{bd}{rK},\tfrac{Kc}{8\pi \xi b})  d\xi
+O(K^{-100})=\frac{1}{c} \sum_{n\ge 1}
n^{-1}W(\tfrac{n}{K},\tfrac{bd}{rK},\tfrac{K^2c}{8\pi
nb})+(K^{-100}),
\end{align}
by Poisson summation again.

As for the outer sum of (\ref{splitresidues}), opening the
Kloosterman sum we have
\begin{align}
&\label{sumram} \sum_{ a \bmod c} S(a^2,b^2,c)e(2ab/c) =\sum_{\substack{a,\gamma \bmod c\\
(\gamma,c)=1}} e\Big(\frac{a^2 \overline{\gamma} + b^2 \gamma + 2ab}{c}\Big)\\
&\nonumber \hspace{1in}=\sum_{\substack{a,\gamma \bmod c\\
(\gamma,c)=1}} e\Big(\frac{\gamma(a+b)^2}{c}\Big) = \sum_{\substack{a,\gamma \bmod c\\
(\gamma,c)=1}} e\Big(\frac{\gamma a^2}{c}\Big)=\sum_{a\bmod c} r_c
(a^2) ,
\end{align}
where $r_c (a^2)$ is a Ramanujan sum and
$\overline{\gamma}\gamma\equiv 1 \bmod c$. We get the last line
above by replacing $a$ by $a\gamma$, then $a$ by $a-b$. Since
$r_c(n) =
\mu\big(\frac{c}{(n,c)}\big)\phi(c)/\phi\big(\frac{c}{(n,c)}\big)$,
we have that the sum of Ramanujan sums in (\ref{sumram}) is
\begin{align}
\phi(c) \sum_{a \bmod c}
\mu\Big(\frac{c}{(a^2,c)}\Big)/\phi\Big(\frac{c}{(a^2,c)}\Big)
\label{sumram2}.
\end{align}
This is multiplicative in $c$, and so we can assume $c$ to be a
prime power $p^k$. Let ord$_p(a)$ denote the highest power of $p$
dividing $a$. If $k$ is even, $\mu\big(\frac{c}{(a^2,c)}\big)=0$
if ord$_p(a) < k/2$, while otherwise
$\mu\big(\frac{c}{(a^2,c)}\big)/\phi\big(\frac{c}{(a^2,c)}\big)=1$.
So if $k$ is even, (\ref{sumram2}) equals $\phi(c)
p^k/p^{k/2}=\phi(c) \sqrt{c}$. For $k$ odd, if ord$_p(a) <
(k-1)/2$ then $\mu\big(\frac{c}{(a^2,c)}\big)=0$, if ord$_p(a) =
(k-1)/2$ then
$\mu\big(\frac{c}{(a^2,c)}\big)/\phi\big(\frac{c}{(a^2,c)}\big)=-\frac{1}{p-1}$
, and if ord$_p(a) > (k-1)/2$ then
$\mu\big(\frac{c}{(a^2,c)}\big)/\phi\big(\frac{c}{(a^2,c)}\big)=1$
. So if $k$ is odd, (\ref{sumram2}) equals
$\phi(c)p^k/p^{(k+1)/2}-\phi(c)
\frac{1}{p-1}(p^k/p^{(k-1)/2}-p^k/p^{(k+1)/2})=0.$
\endproof

\subsection{A Mellin transform}

Denote by $\widetilde{f}$ the Mellin transform of a function $f$.

\begin{lemma} For $0<\Re(s)<1$, we have \label{hmellin}
\begin{align}
\widetilde{\hbar}_{z}(s)= \int_{0}^{\infty}
\frac{h(\sqrt{u})}{\sqrt{2\pi u}}u^{z/2-s}
\Gamma(s)(\cos(s\pi/2)+i\sin(s\pi/2)) du \label{1mellin}
\end{align}
and the bound $\widetilde{\hbar}_{z}(s)\ll_{\Re(z)} (1+|z|)^3
|s|^{-2}.$ For $0<c<1$ we have
\begin{align}
\hbar_z(v)=\frac{1}{2\pi i} \int_{(c)} v^{-s}
\widetilde{\hbar}_{z}(s) ds. \label{2mellin}
\end{align}
\end{lemma}
\proof

We have the Mellin transform $\int_{0}^{\infty} v^{s-1} e^{iv} dv=
\Gamma(s)e^{i\pi s/2}$ for $0<\Re(s)<1$. Thus since $h$ is
compactly supported on $(0,\infty)$ we have
\begin{align}
\int_{0}^{T} v^{s-1} \int_{0}^{\infty}
\frac{h(\sqrt{u})}{\sqrt{2\pi u}} u^{z/2}  e^{iuv} du dv=
\int_{0}^{\infty} \frac{h(\sqrt{u})}{\sqrt{2\pi u}} u^{z/2}
\big(u^{-s}\Gamma(s)e^{i\pi s/2}+o_s(1)\big)du
\end{align}
as $T\rightarrow \infty$. This gives $(\ref{1mellin})$. By
Stirling's approximation we have $\Gamma(s)e^{i\pi s/2}\ll |s|$
for $0<\Re(s)<1$. So on integrating by parts the expression
(\ref{1mellin}) we get
\begin{align}
\widetilde{\hbar}_{z}(s)\ll_{\Re(z)} (1+|z|)^3 |s|^{-2}.
\end{align}
By the bound (\ref{hbarbound}) we have that
$\int_0^{\infty}v^{s-1}\hbar_{z}(v)dv$ converges absolutely, so
that (\ref{2mellin}) follows by Mellin inversion.
\endproof

\subsection{Proof of Theorem \ref{M2}: the diagonal}

Recall the Hecke multiplicative property
\begin{align}
\label{heckerelation} a_f(m^2)a_f(r^2)=\sum_{d|(m^2,r^2)}
a_f(m^2r^2/d^2).
\end{align}
Using the approximate functional equation and the Petersson trace
formula, we get that the twisted second moment is
\begin{align}
&\nonumber \sum_{k \equiv 0 \bmod 2} h\Big(\frac{k-1}{K}\Big)
\summ_{f
\in H_k} L(\hhalf,\sym f)^2 a_f(r^2)\\
&\label{expanded2moment}= 4\sum_{k \equiv 0 \bmod 2}
h\Big(\frac{k-1}{K}\Big) \sum_{m\ge 1}\sum_{d|(r^2,m^2)}\sum_{\substack{n\ge 1\\n=mr/d}} \frac{1}{\sqrt{nm}} V_k(n)V_k(m) \\
&\nonumber + 4 \sum_{k \equiv 0 \bmod 2} h\Big(\frac{k-1}{K}\Big)
\sum_{n,m\ge 1}\sum_{d|(r^2,m^2)}2\pi i^k
\frac{S(n^2,m^2r^2/d^2,c)}{c}J_{k-1}\Big(\frac{4\pi
nmr/d}{c}\Big)\\
&\nonumber \hspace{3in} \times\frac{V_k(n)V_k(m)}{\sqrt{nm}}.
\end{align}
Write $d=d_1d_2^2$, with $d_1$ square-free, and note that $d|m^2$ if
and only if $d_1d_2|m$. Thus we may replace $m$ above by $m d_1 d_2$
and get that (\ref{expanded2moment}) equals
\begin{align}
&\label{expan1}= 4\sum_{k \equiv 0 \bmod 2}
h\Big(\frac{k-1}{K}\Big) \sum_{d|r^2} \sum_{m \ge 1} \frac{1}{m\sqrt{rd_1}} V_k(mr/d_2)V_k(m d _1 d_2) \\
&+ 4 \sum_{k \equiv 0 \bmod 2} h\Big(\frac{k-1}{K}\Big) \sum_{d|r^2}
\sum_{c,n,m\ge 1} 2 \pi i^k
\frac{S(n^2,m^2r^2/d_2^2,c)}{c}J_{k-1}\Big(\frac{4\pi
nmr/d_2}{c}\Big)\label{offdiag}\\
&\nonumber \hspace{2.7in} \times
\frac{V_k(n)V_k(md_1d_2)}{\sqrt{nmd_1d_2}}.
\end{align}

We first evaluate line (\ref{expan1}), the diagonal contribution
to the second moment, and turn to the off-diagonal part
(\ref{offdiag}) in the next section. We have using (\ref{V5}),
\begin{align}
&\nonumber \sum_{m\ge 1} \frac{V_k(mr/d_2)V_k(md_1d_2)}{m}\\
&= \frac{1}{(2\pi i)^2}\int_{(\ep)}\int_{(\ep)}
\frac{G(x)G(y)\zeta(1+2x)\zeta(1+2y)}{xy}\frac{k^{x+y}}{(r/d_2)^{x}
(d_1d_2)^{y}}
 \zeta(1+x+y) dxdy \label{diag} \\
&\nonumber \hspace{4in}+ O(K^{-1+\ep}) .
\end{align}
We move the line of integration $\Re(x)=\ep$ to $\Re(x)=-1+\ep$ ,
crossing a double pole at $x=0$ and a simple pole at $x=-y$. The
integral on the new line is $\ll rK^{-1+\ep}$. Thus since
$\zeta(1+2x)/x=1/(2x^2)+\gamma/x+\ldots$ we get that (\ref{diag}) is
\begin{align}
&\label{maindiag} \frac{1}{2\pi i} \int_{(\ep)}
\frac{G(y)\zeta(1+2y)}{y}\frac{k^y}{(d_1d_2)^y}\Big(\Big(\gamma+
\frac{G^{\prime}(0)}{2}+\frac{\log(kd_2/r)}{2}\Big) \zeta(1+y) +
\frac{\zeta^{\prime}(1+y)}{2}\Big)dy \\
&+\frac{1}{2\pi i} \int_{(\ep)}
\frac{G(-y)G(y)\zeta(1-2y)\zeta(1+2y)}{-y^2}\frac{r^y}{d^y}dy  +
O(rK^{-1+\ep})\label{canceldiag}.
\end{align}
Now we move the line of integration in (\ref{maindiag}) to
$\Re(y)=-1+\ep$, crossing a pole of order 4 at $y=0$. The integral
on the new line is $\ll (K/d_1d_2)^{-1+\ep}.$ From the residue at
the pole we find that (\ref{maindiag}) equals
\begin{align}
&\label{logpolys} \frac{1}{8}\Big(\log
\frac{k}{d_1d_2}\Big)^2\log\frac{kd_2}{r}-\frac{1}{24}\Big(\log
\frac{k}{d_1d_2}\Big)^3
+2P_2\Big(\log\frac{k}{d_1d_2}\Big)+2\log\frac{kd_2}{r}P_1\Big(\log\frac{k}{d_1d_2}\Big)\\
&\nonumber \hspace{4.1in}+O(rK^{-1+\ep}),
\end{align}
where $P_i$ is a polynomial of degree $i$. Inserting
(\ref{logpolys}) into (\ref{expan1}) and averaging over $k$ using
(\ref{hsum}) gives the main term of Theorem \ref{M2}.

Inserting (\ref{canceldiag}) into (\ref{expan1}) gives that the
contribution of this term to the diagonal is
\begin{align}
\frac{-2K}{\sqrt{r}} \hat{h}(0) \frac{1}{2\pi i} \int_{(\ep)}
\frac{G(-y)G(y)\zeta(1-2y)\zeta(1+2y)}{y^2}\sum_{d|r^2}\frac{1}{\sqrt{d_1}}\frac{r^y}{d^y}dy+O(rK^{\ep}).\label{canceldiag2}
\end{align}
We will see that this cancels out with the off-diagonal! Note that
the integrand is an even function of $y$. One can check that
$\sum_{d|r^2}\frac{1}{\sqrt{d_1}}\frac{r^y}{d^y}$ is even by
observing that it is multiplicative in $r$ and then checking that
it is even for $r$ equal to a prime power.

\subsection{Proof of Theorem \ref{M2}: the off-diagonal}

In this section we evaluate line (\ref{offdiag}). By
(\ref{avgwithV}) we have that (\ref{offdiag}) equals
\begin{align}
&-4\pi K \sum_{d|r^2} \sum_{n,m,c\ge 1}
\frac{S(n^2,m^2r^2/d_2^2,c)}{c\sqrt{nmd_1d_2}}\Im\Big(\frac{e(2nmr/cd_2-1/8)}{\sqrt{4\pi
nmr/cd_2}} W_K(n,md_1d_2,\tfrac{K^2cd_2}{8\pi nmr})\Big)\\
&\nonumber \hspace{0.8in} +O\Big(K^{-4}\sum_{d|r^2} \sum_{n,m,c\ge
1} r\sqrt{nm}
\frac{|S(n^2,m^2r^2/d_2^2,c)|}{c^2}\Big(\frac{K}{n}\Big)^{A_1}\Big(\frac{K}{m}\Big)^{A_2}\Big).
\end{align}
The error term above is $O(rK^{-1+\ep})$ by restricting the
summation to $n,m\le K^{1+\ep}$ and using Weil's bound for the
Kloosterman sum. The main term equals by (\ref{wbound}),
(\ref{wtow}) and (\ref{w1bound}),
\begin{align}
&\label{main01}-K\frac{2\sqrt{\pi}}{\sqrt{r}}\Im\Big( \sum_{d|r^2}
\sum_{n,m,c\ge 1} \frac{S(n^2,m^2r^2/d_2^2,c)}{nm\sqrt{cd_1}}
e(2nmr/cd_2-1/8)
W(\tfrac{n}{K},\tfrac{md_1d_2}{K},\tfrac{K^2cd_2}{8\pi
nmr})\Big)\\
&\nonumber \hspace{1.4in}
+O\Big(\frac{K^{\ep}}{\sqrt{r}}\sum_{d|r^2}\sum_{n,m\le
K^{1+\ep}}\sum_{c\le nmr/K^{2-\ep}}
\frac{|S(n^2,m^2r^2/d_2^2,c)|}{nm\sqrt{c}}\Big).
\end{align}
The error above is $O(r^{\half}K^{\ep})$ using Weil's bound on the
Kloosterman sum. Now using Lemma \ref{sumkloo} (with $b=mr/d_2$
and $d=d_1d_2^2$) we get that the main term of (\ref{main01}) is,
up to an error of $O(K^{-10})$,
\begin{align}
\label{offmain1}-K\frac{2\sqrt{\pi}}{\sqrt{r}}\Im\Big(e(-1/8)
\sum_{d|r^2}\frac{1}{\sqrt{d_1}} \sum_{c,n,m\ge 1}
\frac{\phi(c^2)}{c^2}
\frac{1}{nm}W(\tfrac{n}{K},\tfrac{md_1d_2}{K},\tfrac{K^2c^2d_2}{8\pi
nmr})\Big).
\end{align}
The point here is the `completion of the square' which we saw in
(\ref{sumram}). This feature also appears in \cite{iwamich} and
\cite{ils}. Note that the condition $r\le K^{1-\epsilon}$ of Lemma
\ref{sumkloo} is the real limitation to the length of our
mollifier.

We have by Lemma \ref{hmellin} that
\begin{align}
&W(\tfrac{n}{K},\tfrac{md_1d_2}{K},\tfrac{K^2c^2d_2}{8\pi
nmr})=\frac{1}{(2\pi i)^3} \int_{(1)}\int_{(1)} \int_{(1-\ep)}
G(y)G(x) \frac{\zeta(1+2y)}{y}\frac{\zeta(1+2x)}{x}\\
&\nonumber\hspace{2.1in} \times \Big(\frac{K}{n}\Big)^y
\Big(\frac{K}{md_1d_2}\Big)^x \Big(\frac{8\pi
nmr}{K^2c^2d_2}\Big)^s \widetilde{\h}_{x+y}(s) ds dy dx.
\end{align}
Inserting this into (\ref{offmain1}) and exchanging orders of
summation and integration freely by absolute convergence we get that
(\ref{offmain1}) equals
\begin{align}
\label{mellin}&-K\frac{2\sqrt{\pi}}{\sqrt{r}}\Im\Big(e(-1/8)\frac{1}{(2\pi
i)^3} \int_{(1-\ep)}\int_{(1)} \int_{(1)}(8\pi)^sG(y)G(x)
\frac{\zeta(1+2y)}{y}\frac{\zeta(1+2x)}{x}\frac{K^{x+y}}{K^{2s}}\\
&\nonumber\hspace{0.4in}
\times\zeta(1+x-s)\zeta(1+y-s)\zeta(2s)\zeta(2s+1)^{-1}\widetilde{\h}_{x+y}(s)\sum_{d|r^2}\frac{r^s}{d_1^xd_2^{x+s}\sqrt{d_1}}dy
dx ds\Big),
\end{align}
where we used that $\zeta(2s)\zeta(2s+1)^{-1}=\sum_{c\ge 1}
\frac{\phi(c^2)}{c^{2s+2}}$ for $\Re(s)>\half$. Now we move the
lines of integration $\Re (x) = \Re(y) =1$ to $\Re(x)=\Re(y)= \ep$,
crossing simple poles at $y=s$ and $x=s$. On the new lines
$|K^{x-s}|\ll K^{-1+\epsilon}$, $|K^{y-s}|\ll K^{-1+\epsilon}$, so
the contribution to (\ref{mellin}) of the new integrals is
$O(r^{\half}K^{\epsilon}).$ Thus from the poles we get that
(\ref{mellin}) equals
\begin{align}
\label{aline}&-K\frac{2\sqrt{\pi}}{\sqrt{r}}\Im\Big(e(-1/8)\frac{1}{2\pi
i} \int_{(1-\ep)}(8\pi)^sG(s)^2 \frac{\zeta(1+2s)^2}{s^2}
\zeta(2s)\zeta(2s+1)^{-1}\widetilde{\h}_{2s}(s)\\
&\nonumber\hspace{2.4in}
\times\sum_{d|r^2}\frac{r^s}{d_1^sd_2^{2s}\sqrt{d_1}}ds\Big)
+O(r^{\half}K^{\ep}).
\end{align}
Using Lemma \ref{hmellin} the main term of (\ref{aline}) equals
\begin{align}
&\nonumber-K\frac{2\sqrt{\pi}}{\sqrt{r}}\Im\Big(e(-1/8)\frac{1}{2\pi
i} \int_{0}^{\infty}\int_{(1-\ep)}
\frac{(8\pi)^sG(s)^2\zeta(1+2s)\zeta(2s)}{s^2}\frac{h(\sqrt{u})}{\sqrt{2\pi
u}} \Gamma(s)
\\
&\nonumber\hspace{2.4in}\times\big(\cos(\tfrac{\pi
s}{2})+i\sin(\tfrac{\pi s}{2})\big)
\sum_{d|r^2}\frac{1}{\sqrt{d_1}}\frac{r^s}{d^s} ds du\Big)\\
&\label{cancels}=-K\frac{2\sqrt{\pi}}{\sqrt{r}} \int_{0}^{\infty}
\frac{h(\sqrt{u})}{\sqrt{2\pi u}} du \frac{1}{2\pi
i}\int_{(1-\ep)}
\frac{(8\pi)^sG(s)^2\zeta(1+2s)\zeta(2s)}{s^2}\Gamma(s)\sin(\tfrac{\pi
s}{2}-\tfrac{\pi}{4})\\
&\nonumber\hspace{2.9in} \times
\sum_{d|r^2}\frac{1}{\sqrt{d_1}}\frac{r^s}{d^s} ds.
\end{align}
We will see by the functional equation of $\zeta(s)$ that the
integrand above is an even function of $s$. In fact
(\ref{cancels}) exactly cancels out with (\ref{canceldiag2}). This
is another feature of the moment calculation.

We have $-2\sqrt{\pi}\int_{0}^{\infty} \frac{h(\sqrt{u})}{\sqrt{2\pi
u}} du=-2\sqrt{2}\hat{h}(0)$, so to show that (\ref{cancels}) and
(\ref{canceldiag2}) cancel we need to show that
$(8\pi)^s\zeta(2s)G(s)\Gamma(s)\sin(\tfrac{\pi
s}{2}-\tfrac{\pi}{4})=-\frac{1}{\sqrt{2}}\zeta(1-2s)G(-s)$. By the
functional equation of $\zeta(s)$ and the duplication formula of
$\Gamma(s)$ we have
\begin{align}
\label{combine1}
\Gamma(s)\zeta(2s)=\pi^{2s-\half}\Gamma(\hhalf-s)\zeta(1-2s)=2^{-\half-s}\pi^{2s-1}\Gamma(\tfrac{1}{4}-\tfrac{s}{2})\Gamma(\tfrac{3}{4}-\tfrac{s}{2})\zeta(1-2s).
\end{align}
By the definition of $G(s)$ given in Lemma \ref{approxeq} and the
reflection formula of $\Gamma(s)$ we have
\begin{align}
\label{combine2} G(s)\sin(\tfrac{\pi
s}{2}-\tfrac{\pi}{4})=\frac{\pi^{-\frac{3}{2}s}2^{-s}\Gamma(\frac{s}{2}+\frac{3}{4})\sin(\tfrac{\pi
s}{2}-\tfrac{\pi}{4})}{\Gamma(\frac{3}{4})}=\frac{-\pi^{1-\frac{3}{2}s}2^{-s}}{\Gamma(\frac{3}{4})\Gamma(\frac{1}{4}-\frac{s}{2})}.
\end{align}
Combining (\ref{combine1}) and (\ref{combine2}) gives the required
formula for $(8\pi)^s\zeta(2s)G(s)\Gamma(s)\sin(\tfrac{\pi
s}{2}\break-\tfrac{\pi}{4})$.

\section{Optimizing the mollifier}

In the previous sections we found the first and second moments of
$L(\half,\sym f)$ twisted by the Fourier coefficient $a_f(r^2)$.
These immediately lead to the moments of $L(\half,\sym f)$ times
the mollifier $M(f)=\sum_{r} a_f(r^2) x_r r^{-\half}$. Here $x_r$
are coefficients to be determined with $x_r=0$ for $r\le 0$ or $r>
M$, where $M=K^a$ for some constant $a>0$ also to be determined.
Also suppose that $x_r \ll M^{\ep}$ as this will be the case once
choices are made. In this chapter we choose a mollifier that
maximizes the proportion of non-vanishing of $L(\half,\sym f)$
given by (\ref{ratiointro}). We shall see that this boils down to
a problem of minimizing a quadratic form in the $x_r$ with a
linear constraint. By Theorem \ref{M1} we have the first mollified
moment, a linear form in $x_r$:
\begin{align}
{\mathcal M}_1&=\sum_{k \equiv 0 \bmod 2} h\big(\frac{k-1}{K}\big)
\summ_{f \in H_k} L(\hhalf,\sym f) \sum_{r} \frac{x_r}{\sqrt{r}}a_f(r^2)\\
&\nonumber= \frac{K}{2}\sum_{r} \frac{x_r}{r} \int_{0}^{\infty}
h(u)(\log (uK/r)+C)du+O(MK^\ep).
\end{align}
Using the Hecke relation (\ref{heckerelation}) we have the second
mollified moment, a quadratic form in $x_r$:
\begin{align}
{\mathcal M}_2&=\sum_{k \equiv 0 \bmod 2} h\big(\frac{k-1}{K}\big)
\summ_{f \in H_k} L(\hhalf,\sym f)^2 \Big(\sum_{r} \frac{x_r}{\sqrt{r}}a_f(r^2)\Big)^2\\
&\nonumber=\sum_{k \equiv 0 \bmod 2} h\big(\frac{k-1}{K}\big)
\summ_{f \in H_k} L(\hhalf,\sym f)^2 \sum_{\al\ge1}\sum_{e|\al^2}
\sum_{\substack{r_1,r_2\\(r_1,r_2)=1}}\frac{x_{\al r_1}x_{\al
r_2}}{\al \sqrt{r_1 r_2}}
a_f\Big(\frac{r_1^2r_2^2\al^4}{e^2}\Big).
\end{align}
Above $\frac{r_1 r_2\al^2}{e}\le K^{1-\ep}$ if $M\le K^{\half-\ep}$.
Thus under this condition we have by Theorem \ref{M2},
\begin{align}
\label{calm2}{\mathcal M}_2=\frac{K}{2}\sum_{\al}\sum_{e|\al^2}
\sum_{\substack{r_1,r_2\\(r_1,r_2)=1}}x_{\al r_1} x_{\al
r_2}\frac{\sqrt{e}}{\al^2 r_1 r_2} S\Big(\frac{r_1 r_2
\al^2}{e}\Big)+O(M^2K^{\ep}),
\end{align}
where
\begin{align} \label{slog} S(r)=\sum_{d|r^2} \frac{1}{\sqrt
d_1} \int_{0}^{\infty} &h(u)
\Big(\frac{1}{2}\log^2\frac{uK}{d_1d_2}\log\frac{uKd_2}{r}-\frac{1}{6}\log^3\frac{uK}{d_1d_2}
+2 P_2\Big(\log\frac{uK}{d_1d_2}\Big)\\
&\nonumber\indent\indent\indent+2
\log\frac{uKd_2}{r}P_1\Big(\log\frac{uK}{d_1d_2}\Big) \Big) du,
\end{align}
and $P_i$ are degree $i$ polynomials as in Theorem \ref{M2}. We
see from this that the longest length of mollifier we are able to
take is $M=K^{\half-\ep}$ and hence we make the restriction
$a<\half$. Define a multiplicative function
\begin{align}
\tau(n)=\sum_{d|n^2} \frac{1}{\sqrt d_1}=\sum_{f|n}
\frac{\mu^2(f)}{\sqrt{f}} d(n/f),
\end{align}
where in the second sum $d$ is the divisor function. Note that
$\tau(p)=2+p^{-\half}$, so that $\tau$ behaves roughly like the
divisor function. Since $S(r)$ is roughly of size $(\log K)^3$ for
$r\le K^{1-\ep}$, we expect ${\mathcal M}_2$ to be roughly of size
\begin{align}
\label{m2quad}K (\log K)^3 \sum_{\al}\sum_{e|\al^2}
\sum_{\substack{r_1,r_2\\(r_1,r_2)=1}}\frac{\sqrt{e}}{\al^2 r_1
r_2} \tau\Big(\frac{r_1 r_2 \al^2}{e}\Big) x_{\al r_1}x_{\al r_2}
\end{align}
We will make a choice of $x_r$ based on this quadratic form; one
that maximizes the ratio of ${\mathcal M}_1^2$ to (\ref{m2quad}),
under the assumption that $x_r$ is supported on square-free $r$.
The form of our mollifier and this last assumption are reasonable
from the Euler product $(\ref{euler})$ and the expectation that
the mollifier will mimic a multiple of the inverse of
$L(\half,\sym f)$. For this choice we will find that ${\mathcal
M}_1 \asymp K(\log K)^3$ and ${\mathcal M}_2 \asymp K(\log K)^6$,
and we'll get that
\begin{align}
\frac{{\mathcal M}_1^2}{{\mathcal M}_2}\sim
\Big(1-\frac{1}{(1+a)^{3}}\Big)\sum_{k \equiv 0 \bmod 2}
h\Big(\frac{k}{K}\Big)\summ_{f\in H_k} 1.
\end{align}
Thus together with observation (\ref{ratiointro}) the main theorem
will follow.

\subsection{Change of variables}

We simplify the problem by diagonalizing the quadratic form
(\ref{m2quad}). We assume that $x_r=0$ if $r$ is not square-free.
We have
\begin{align}
&\nonumber\sum_{\al}\sum_{e|\al^2}
\sum_{\substack{r_1,r_2\\(r_1,r_2)=1}}\frac{\sqrt{e}}{\al^2 r_1 r_2}
\tau\Big(\frac{r_1 r_2 \al^2}{e}\Big) x_{\al r_1}x_{\al r_2}\\
&\nonumber=\sum_{\al}\sum_{e|\al^2}
\sum_{\substack{r_1,r_2\\(r_1,r_2)=1}}\frac{\sqrt{e}}{\al^{2}r_1r_2}
\tau(r_1)\tau(r_2)\tau(\al^2/e)x_{\al
r_1}x_{\al r_2}\\
&\label{quadnondiag}=\sum_{\al,\beta \ge1}\sum_{e|\al^2}
\sum_{r_1,r_2}\frac{\mu(\beta)\sqrt{e}}{\beta^2\al^{2}r_1r_2}
 \tau(\beta)^2
\tau(r_1)\tau(r_2)\tau(\al^2/e)x_{\al \beta r_1}x_{\al \beta r_2},
\end{align}
where in the last line we used the Mobius function to detect the
condition $(r_1,r_2)=1$. Define a change of variables
\begin{align}
y_j= \sum_r \frac{\tau(r)}{r} x_{jr}.
\end{align}
Thus $y_j$ is supported on positive square-free integers $j$ less
than or equal to $M$. This change of variables is invertible:
\begin{align}
x_r=\sum_n
\frac{x_{nr}\tau(n)}{n}\sum_{d|n}\mu(d)=\sum_{d}\frac{\tau(d)\mu(d)}{d}\sum_{n}
\frac{x_{ndr}\tau(n)}{n}=\sum_d \frac{y_{dr}\tau(d)\mu(d)}{d},
\end{align}
where we used the Mobius function to detect $n=1$ in the first
sum. Also define coefficients for positive square-free $j$ (and
otherwise set equal to zero),
\begin{align}
v_j=\frac{1}{j^{2}} \sum_{\al \beta=j} \mu(\beta) \tau(\beta)^2
\sum_{e|\al^2} \sqrt{e}
\tau\Big(\frac{\al^2}{e}\Big)=\frac{1}{j}\prod_{p|j}\Big(1+\frac{2}{p^{\half}}-\frac{2}{p^{\frac{3}{2}}}-\frac{1}{p^2}\Big).
\end{align}
Then we have that (\ref{quadnondiag}) equals
\begin{align}
\sum_{j} v_j y_j^2,
\end{align}
and (\ref{m2quad}) equals $K\log^3 K \sum_{j} v_j y_j^2$. This is
the diagonal quadratic form that we will work with. For this change
of variables, the linear form ${\mathcal M}_1$ equals
\begin{align}
{\mathcal M}_1=\frac{\hat{h}(0)K}{2}\sum_j u_j y_j+O(M^2K^{\ep}),
\end{align}
where the coefficients are
\begin{align}
u_j=\frac{1}{\hat{h}(0)}\frac{1}{j}\sum_{nr=j} \mu(n) \tau(n)
\int_{0}^{\infty} h(u)(\log(uK/r)+C) du.
\end{align}
For square-free $j$ we have
\begin{align}
u_j=\frac{1}{j} \sum_{n|j} \mu(n) \tau(n) \log(Kn/j)
+O(j^{-1}\mu_*(j)),
\end{align}
where we define for positive square-free $j$ (and otherwise set
equal to zero),
\begin{align}
\mu_*(j)=\sum_{n|j} \mu(n)
\tau(n)=\prod_{p|j}\Big(-1-\frac{1}{\sqrt{p}}\Big).
\end{align}
As the notation suggests, $\mu_*$ behaves roughly like $\mu$. Now,
\begin{align}
&\nonumber\sum_{n|j} \mu(n) \tau(n) \log(Kn/j)=\mu_*(j)\log
K-\sum_{n|j} \mu(n)
\tau(n) \sum_{p|\frac{j}{n}} \log p\\
&\nonumber=\mu_*(j)\log K-\sum_{p|j} \log p\sum_{n|\frac{j}{p}}
\mu(n)
\tau(n)=\mu_*(j)\log K-\mu_*(j) \sum_{p|j}  \frac{\log p}{\mu_*(p)}\\
&\nonumber=\mu_*(j)\log K-\mu_*(j) \sum_{p|j} \log p(-1+O(p^{-\half}))\\
&=\label{ujval} \mu_*(j)\log
(Kj)+O\Big(\prod_{p|j}(1+p^{-\frac{1}{2}+\ep})\Big) .
\end{align}
Thus
\begin{align}
u_j=\frac{\mu_*(j)}{j}\log(Kj)+O\Big(j^{-1}\prod_{p|j}(1+p^{-\frac{1}{2}+\ep})\Big).
\end{align}
With this change of variables we can easily determine the optimum
choice for $y_j$. By the Cauchy-Schwarz inequality we have,
\begin{align}
\Big(\sum_j u_j y_j\Big)^2=\Big(\sum_j \frac{u_j}{\sqrt{v_j}}
\sqrt{v_j} y_j\Big)^2 \le \Big(\sum_j v_j y_j^2\Big) \Big(\sum_j
\frac{u_j^2}{v_j}\Big).
\end{align}
Thus $(\sum_j u_j y_j)^2/\sum_j v_j y_j^2$ obtains its upper bound
when
\begin{align}
y_j=\frac{u_j}{v_j}.\label{choice}
\end{align}
Set $y_j$ to this value henceforth and note that $y_j\ll \log K.$

\subsection{Evaluating ${\mathcal M}_1$}

For the choice made in (\ref{choice}) we have
\begin{align}
\nonumber{\mathcal M}_1&= \frac{\hat{h}(0)K}{2} \sum_{j\le M} \frac{u_j^2}{v_j}\\
&\label{m1error}=\frac{\hat{h}(0)K}{2}\sum_{j\le M}
\frac{\mu_*(j)^2}{j^2v_j}\log^2(Kj)+O\Big(K\log K\sum_{j\le M}
j^{-1}\prod_{p|j}(1+p^{-\frac{1}{2}+\ep})\Big).
\end{align}
The error is $O(K\log ^2 K)$. We have that
\begin{align}
\label{easyfactor}
\frac{\mu_*(p)^2}{p^2v_p}=\frac{1+2p^{-\half}+p^{-1}}{p(1+2p^{-\half}-2p^{-\frac{3}{2}}-p^{-2})}=\frac{1}{p(1-p^{-1})}=\frac{1}{\phi(p)},
\end{align}
after a nice simplification. Now it is well known that
\begin{align}
\label{asympuv} \sum_{j\le x} \frac{\mu_*(j)^2}{j^2v_j}=
\sum_{j\le x} \frac{\mu(j)^2}{\phi(j)} =\log x + O(1).
\end{align}
Indeed, we have that $\sum_{j\ge 1} \frac{\mu(j)^2}{\phi(j)}
\frac{1}{j^s}=\prod_p
(1+p^{-(1+s)}+p^{-(2+s)}+\ldots)=\zeta(s+1)F(s)$, where $F(s)$ is
analytic and absolutely convergent for $\Re(s)>-\half$ and
$F(0)=1$.

By partial summation we get
\begin{align}
\nonumber{\mathcal M}_1&=\frac{\hat{h}(0)K}{2}\sum_{j\le M}
\frac{\mu_*(j)^2}{j^2v_j}\log^2(Kj)+O(K \log^2 K)\\
&\nonumber= \frac{\hat{h}(0)K}{2}\log^2(KM)\log M-\frac{\hat{h}(0)K}{2} \sum_{j\le M}(\log^2(K(j+1))- \log^2(Kj))\log j\\
&\nonumber\hspace{0.5in}+O(K\log^2 K)\\
&\sim \hat{h}(0)K\log^3 K \frac{a(3+3a+a^2)}{6} .
\end{align}

\subsection{Evaluating ${\mathcal M}_2$} Expanding the logarithms in
the definition of $S(r_1r_2\al^2/e)$, we see that ${\mathcal M}_2$
is a sum of terms proportional to
\begin{align}
\label{otherterms} &K(\log K)^{i_0} \sum_{\al}\sum_{e|\al^2}
\sum_{\substack{r_1,r_2\\(r_1,r_2)=1}}\frac{\sqrt{e}}{ r_1 r_2}
\frac{\log(\al^2/e)^{i_1}}{\al^2} \sum_{d|r_1^2r_2^2\al^4/e^2}
\frac{(\log d_1)^{i_2}(\log
d_2)^{i_3}}{\sqrt{d_1}}\\
&\nonumber\hspace{1.6in} (\log r_1 r_2)^{i_4} \Big(
\int_{0}^{\infty}h(u)(\log u)^{i_5} du\Big) x_{\al r_1}x_{\al
r_2},
\end{align}
where $i_0,\ldots, i_5$ are non-negative integers with
$i_0+\ldots+i_5= 3$, $i_1\le 1$ and $i_4\le 1$. We show in this
section that (\ref{otherterms}) is of size $K(\log K)^6$ when
$i_0+i_3+i_4=3$. In the next section we will see that it is
bounded by $K(\log K)^5$ otherwise. Thus in this section we want
to evaluate the quadratic form,
\begin{align}
\label{mainquad} \frac{K\hat{h}(0)}{2}\sum_{\al}\sum_{e|\al^2}
\sum_{\substack{r_1,r_2\\(r_1,r_2)=1}}\frac{\sqrt{e}}{\al^2 r_1
r_2}
\sum_{d|r_1^2r_2^2\al^4/e^2}\Big(\half\log^2\frac{K}{d_2}\log\frac{Kd_2}{r_1r_2}-\frac{1}{6}\log^3\frac{K}{d_2}\Big)x_{\al
r_1}x_{\al r_2},
\end{align}
for the choice made in (\ref{choice}). We have
\begin{align}
\label{logterms}
\half\log^2\frac{K}{d_2}\log\frac{Kd_2}{{r_1r_2}}-\frac{1}{6}\log^3\frac{K}{d_2}=&\frac{1}{3}\log^3
K-\half \log^2 K \log r_1r_2+\log K \log r_1r_2 \log d_2\\
&\nonumber\indent - \log K \log^2 d_2 + \frac{2}{3} \log^3 d_2
-\half \log r_1r_2 \log^2 d_2,
\end{align}
and we evaluate the contribution of each of these terms one by
one.

$\bullet$ The contribution of the $\log^3 K$ term of
(\ref{logterms}) to (\ref{mainquad}) is
\begin{align}
\label{r1} \frac{K\hat{h}(0)}{6}\sum_{j} v_j
y_j^2=\frac{K\hat{h}(0)}{6}\sum_{j\le M} \frac{u_j^2}{v_j}\sim
\hat{h}(0)K\log^6 K\frac{a(3+3a+a^2)}{18},
\end{align}
by the same calculations as for the evaluation of ${\mathcal M}_1$.

$\bullet$ The contribution of the $\log^2 K \log r_1r_2$ term of
(\ref{logterms}) to (\ref{mainquad}) is
\begin{align}
&\nonumber \frac{-\hat{h}(0)K\log^2 K}{4} \sum_{\al}\sum_{e|\al^2}
\sum_{\substack{r_1,r_2\\(r_1,r_2)=1}}\frac{\sqrt{e}}{\al^2 r_1
r_2} \tau\Big(\frac{r_1 r_2 \al^2}{e}\Big) \log(r_1r_2) x_{\al
r_1}x_{\al
r_2}\\
&\label{log2klogr00} =\frac{-\hat{h}(0)K\log^2 K}{2}
\sum_{\al}\sum_{e|\al^2}
\sum_{\substack{r_1,r_2\\(r_1,r_2)=1}}\frac{\sqrt{e}}{\al^2 r_1
r_2} \tau\Big(\frac{r_1 r_2 \al^2}{e}\Big) \log(r_1) x_{\al
r_1}x_{\al r_2},
\end{align}
by symmetry. Using the Mobius function to detect $(r_1,r_2)=1$ this
equals
\begin{align}
&\nonumber \frac{-\hat{h}(0)K\log^2 K}{2} \sum_{\al,\beta}
\sum_{e|\al^2} \frac{\mu(\beta)
\tau(\beta)^2\sqrt{e}\tau(\al^2/e)}{\al^2\beta^2}\sum_{r_1,r_2}
\frac{\tau(r_1)\tau(r_2) \log(\beta r_1)}{r_1r_2} x_{\al \beta
r_1}x_{\al \beta
r_2}\\
&\label{log2klogr}=\frac{-\hat{h}(0)K\log^2 K}{2}\sum_j  v_j y_j
\sum_{r_1}
\frac{\tau(r_1)\log(r_1)}{r_1}x_{j r_1}\\
&\nonumber\hspace{0.8in}+ \frac{-\hat{h}(0)K\log^2 K}{2}
\sum_{\al,\beta} \Big(\frac{1}{(\al \beta)^2} \mu(\beta)
\tau(\beta)^2 \log \beta \sum_{e|\al^2} \sqrt{e} \tau(\al^2/e)
\Big) y_{\al \beta}^2,
\end{align}
in terms of the new variables $y_j$. Using that $\al^{-2}
\sum_{e|\al^2} \sqrt{e} \tau(\al^2/e)\ll
\al^{-1}\prod(1+p^{-\frac{1}{2}+\ep})$ and the bound $y_j\ll \log
K$, we have that the second line of (\ref{log2klogr}) is $O(\log^5
K)$.

Let $\Lambda_j(n)=\sum_{d|n}\mu(d)(\log\tfrac{n}{d})^j$ for $j\ge
0$ be the generalized von Mangoldt function. Thus $\Lambda_0(1)=1$
and $\Lambda_0(n)=0$ for $n\ge 2$, and $\Lambda_1(n)$ is the usual
von Mangoldt function $\Lambda(n)$ supported on prime powers. The
function $\Lambda_j(n)$ is supported on integers $n$ having at
most $j$ distinct prime factors and satisfies
$\Lambda_j(n)\ll_j(\log n)^j$, $\sum_{n\le x} \Lambda_j(n)/n \ll_j
(\log  x)^j$ and $(\log n)^j=\sum_{d|n} \Lambda_j(d)$. We also
have $\Lambda_j(nm)=\sum_{0\le i\le j} {j \choose
i}\Lambda_i(n)\Lambda_{j-i}(m)$ for coprime integers $n$ and $m$.
Writing $\log r_1= \sum_{b|r_1} \Lambda(b)$, the first term of
(\ref{log2klogr}) is
\begin{align}
&\nonumber \frac{-\hat{h}(0)K\log^2 K}{2}\sum_j  v_j \sum_b
\frac{\Lambda(b)\tau(b)}{b} y_{bj}y_j\\
&=\label{log2klogrfin}\frac{-\hat{h}(0)K\log^2 K}{2}\sum_{j\le M}
\frac{\mu_*(j)^2}{j^2v_j} \sum_{\substack{b\le M/j\\ (b,j)=1}}
\frac{\Lambda(b)\tau(b)}{b} \frac{\mu_*(b)}{b
v_b}\log(Kj)\log(Kjb)\\
&\nonumber\indent\indent+O(K\log^{5}K),
\end{align}
on substituting the value $u_j/v_j$ of $y_j$ and using
(\ref{ujval}). As the sum above is restricted to prime values of
$b$, we can remove the condition $(b,j)=1$ since otherwise $b|j$
and the contribution of such terms is $\ll K \log^{5}K$. The above
sum can be evaluated by partial summation. Recall that
$\frac{\mu_*(j)^2}{j^2v_j}=\frac{\mu(j)^2}{\phi(j)}$ behaves like
$\frac{1}{j}$ on average, and for prime $b$,
$\tau(b)=2+b^{-\frac{1}{2}}$, $\mu_*(b)=-1-b^{-\frac{1}{2}}$ and
$bv_b=1+O(b^{-\frac{1}{2}})$. Using the prime number theorem we
get that (\ref{log2klogrfin}) is
\begin{align}
\nonumber &\sim \frac{-\hat{h}(0)K\log^2 K}{2} \sum_{j\le M}
\frac{\mu_*(j)^2}{j^2v_j}\Big(-2\log(M/j) \log^2(Kj) -\log^2(M/j)
\log(Kj)\Big)\\
&\label{r2} \sim \hat{h}(0)K\log^6 K \frac{a^2(2+a)^2}{8},
\end{align}
by partial summation.

$\bullet$ The contribution of the $\log K \log r_1r_2 \log d_2$
term of (\ref{logterms}) to (\ref{mainquad}) is
\begin{align}
&\label{logklogrlogd1} \hat{h}(0)K\log K \sum_{\al}\sum_{e|\al^2}
\sum_{\substack{r_1,r_2\\(r_1,r_2)=1}}\frac{\sqrt{e}}{\al^2 r_1 r_2}
\Big(\sum_{d|r_1^2r_2^2\al^4/e^2}\frac{\log d_2}{\sqrt{d_1}}  \Big)
\log r_1 x_{\al r_1}x_{\al r_2}.
\end{align}
Note that
\begin{align}
\label{use}\sum_{d|n^2}\frac{\log^{i}d_1\log^{j}d_2}{\sqrt{d_1}}=\sum_{f|n}
\frac{\mu(f)^2 \log^{i}f}{\sqrt{f}}\sum_{e|\frac{n}{f}} \log^{j}e
=\sum_{a|n}\sum_{b|\frac{n}{a}}\frac{\mu(a)^2
\Lambda_{i}(a)\Lambda_{j}(b)}{\sqrt{a}}\tau\Big(\frac{n}{ab}\Big).
\end{align}
Using this we get
\begin{align}
 \sum_{d|r_1^2r_2^2\al^4/e^2}\frac{\log
d_2}{\sqrt{d_1}}=\Big(\sum_{c|r_1}&\Lambda(c)\tau(r_1/c)\Big)\tau(r_2)\tau(\al^2/e)+\tau(r_1)\Big(\sum_{c|r_2}\Lambda(c)\tau(r_2/c)\Big)\tau(\al^2/e)
\\
&
\nonumber+\tau(r_1)\tau(r_2)\Big(\sum_{c|\al^2/e}\Lambda(c)\tau(\al^2/ec)\Big).
\end{align}
Thus (\ref{logklogrlogd1}) is
\begin{align}
&\label{logklogrlogd1mid}\hat{h}(0)K\log
K\\
&\nonumber \times \Big(\sum_{\al,\beta}\sum_{e|\al^2}\sum_c
\frac{\sqrt{e}\tau(\beta)^2\mu(\beta)\tau(\al^2/e)}{\al^2\beta^{2}}\frac{\Lambda(c)}{c}\sum_{\substack{r_1,r_2\\c\nmid
 r_2}}\frac{x_{\al\beta c r_1} x_{\al\beta r_2}
\tau(r_1)\tau(r_2)\log(c\beta r_1)}{r_1
r_2}\\
&\nonumber+\sum_{\al,\beta}\sum_{e|\al^2}\sum_c
\frac{\sqrt{e}\tau(\beta)^2\mu(\beta)\tau(\al^2/e)}{\al^2\beta^{2}}\frac{\Lambda(c)}{c}\sum_{\substack{r_1,r_2\\c\nmid
r_1}}\frac{x_{\al\beta r_1} x_{\al\beta c r_2}\tau(r_1)\tau(r_2)
\log(\beta r_1)}{r_1
r_2}\\
&\nonumber+\sum_{\al,\beta}\sum_{e|\al^2} \Big(
\frac{\sqrt{e}\tau(\beta)^2\mu(\beta)}{\al^2\beta^{2}}\sum_{c|\al^2/e}\Lambda(c)\tau(\al^2/ec)\Big)\sum_{r_1,r_2}\frac{x_{\al\beta
r_1} x_{\al\beta r_2} \tau(r_1)\tau(r_2)\log(\beta r_1)}{r_1
r_2}\Big).
\end{align}
In terms of the new variables $y_j$, the first two lines of
(\ref{logklogrlogd1mid}) contribute
\begin{align}
&\nonumber\hat{h}(0)K\log K\sum_j v_j \Big(\sum_{c}
\frac{\Lambda(c)\log c}{c} y_{cj}y_{j} +\sum_{c,b}
\frac{\Lambda(c)}{c}
\frac{\Lambda(b)\tau(b)}{b} y_{cbj}y_{j} \\
&\nonumber\hspace{1.5in}+\sum_{c,b}  \frac{\Lambda(c)}{c}
\frac{\Lambda(b)\tau(b)}{b}
y_{cj}y_{bj}\Big)+O(K\log^5 K)\\
 &\sim \hat{h}(0)K\log
K\sum_{j\le M}\frac{\mu_*(j)^2}{j^2 v_j}\Big(\sum_{c\le M/j}
\frac{\Lambda^2(c)}{c}\frac{\mu_*(c)}{cv_c}\log(Kjc)\log(Kj)\\
&\nonumber\hspace{0.8in}+\sum_{bc\le M/j}
\frac{\Lambda(c)\Lambda(b)\tau(b)}{cb} \frac{\mu_*(cb)}{cbv_{cb}}  \log(Kjcb)\log(Kj)\\
&\nonumber\hspace{1.6in}+\sum_{\substack{c\le M/j\\b\le M/j}}
\frac{\Lambda(c)\Lambda(b)\tau(b)}{cb}
\frac{\mu_*(c)\mu_*(b)}{cv_cbv_b} \log(Kjc)\log(Kjb)\Big).
\end{align}
By the prime number theorem this is
\begin{align}
\label{logklogrlogdfinal}&\sim \hat{h}(0)K\log K\sum_{j\le
M}\frac{\mu_*(j)^2}{j^2 v_j}\Big(-\half
\log^2(M/j)\log^2(Kj)-\frac{1}{3}\log^3(M/j)\log(Kj)\\
&\nonumber\hspace{1.8in}+2\Big(\log(M/j)\log(Kj)+\half\log^2(M/j)\Big)^2\\
&\nonumber\hspace{2in}+\log^2(M/j)\log^2(Kj)+\frac{2}{3}\log^3(M/j)\log(Kj)\Big).
\end{align}
We will evaluate (\ref{logklogrlogdfinal}) by partial summation
after combining it with (\ref{logklog2dfinal}), the contribution of
the $\log K \log^2 d_2$ term. As for the third line of
(\ref{logklogrlogd1mid}), we have
\begin{align}
&\sum_{e|\al^2}\frac{\sqrt{e}}{\al^2}\sum_{c|\al^2/e}\Lambda(c)\tau(\al^2/ec)=\frac{1}{\al}\sum_{c|\al^2}\frac{\Lambda(c)}{\sqrt{c}}\sum_{e|\al^2/c}
\frac{\tau(e)}{\sqrt{e}}\ll
\frac{1}{\al}\prod_{p|\al}(1+p^{-\frac{1}{2}+\ep}).
\end{align}
So the contribution of the third line is $O(K \log^5 K)$.

$\bullet$ The contribution of the $\log K \log^2 d_2$ term of
(\ref{logterms}) to (\ref{mainquad}) is
\begin{align}
&\label{logklog2d}\frac{-\hat{h}(0)K\log K}{2}
\sum_{\al}\sum_{e|\al^2}
\sum_{\substack{r_1,r_2\\(r_1,r_2)=1}}\frac{\sqrt{e}}{\al^2 r_1 r_2}
\Big(\sum_{d|r_1^2r_2^2\al^4/e^2}\frac{\log^2 d_2}{\sqrt{d_1}} \Big)
x_{\al r_1}x_{\al r_2}.
\end{align}
We have
\begin{align}
&\sum_{d|r_1^2r_2^2\al^4/e^2}\frac{\log^2
d_2}{\sqrt{d_1}}\\
&\nonumber\indent=\sum_{\substack{i+j+k=2\\i,j,k\ge
0}}a_{ijk}\Big(\sum_{c|r_1}\Lambda_i(c)\tau(r_1/c)\Big)\Big(\sum_{c|r_2}\Lambda_j(c)\tau(r_2/c)\Big)\Big(\sum_{c|\al^2/e}\Lambda_k(c)\tau(\al^2/ec)\Big),
\end{align}
for some coefficients $a_{ijk}$, with $a_{110}=2$. Only the term
with $i=j=1$ and $k=0$ of this expansion is significant. Showing
that the contribution to (\ref{logklog2d}) of any term with $k\ge
1$ is $O(K\log^{5}K)$ is very similar to how we showed that third
term of (\ref{logklogrlogd1mid}) is small. This leaves us to
consider the contribution of the term with $i=2$ and $j=k=0$. This
is less than
\begin{align}
& K\log K\sum_j v_j \sum_c
\frac{\Lambda_2(c)}{c}y_{cj}y_j+O(K\log^{5}K)\\
&\ll \nonumber K\log K \sum_{j\le M}
\frac{\mu_*(j)^2}{jv_j}\sum_{c\le M/j}
\frac{\Lambda_2(c)}{c}\frac{\mu_*(c)}{cv_c} \log(Kjc)\log(Kj)+O(K\log^{5}K)\\
&\nonumber \ll K\log^{5}K,
\end{align}
where the last step follows by the estimate $\sum_{c\le
x}\frac{\Lambda_2(c)\mu_*(c)}{c^2v_c}\ll \log x$. To see this
estimate, note that for prime values of $c$ the sum is $\sum_{p\le
x}\frac{-(\log p)^2}{p} (1+O(\frac{1}{\sqrt{p}}
))=-\frac{1}{2}(\log x)^2+ O(\log x)$ by the prime number theorem,
and for $c$ a product of two distinct primes the sum is
$\sum_{\substack{pq\le x\\p\neq q}}\frac{2 \log p \log q
}{pq}(1+O(\frac{1}{\sqrt{p}}+\frac{1}{\sqrt{q}}
))=\frac{1}{2}(\log x)^2+ O(\log x).$

Thus (\ref{logklog2d}) is
\begin{align}
\nonumber &\sim \frac{-\hat{h}(0)K\log K}{2} \sum_j 2v_j \Big(\sum_c
\frac{\Lambda(c)}{c} y_{jc}\Big)^2\\
&\nonumber\sim -\hat{h}(0)K\log K \sum_{j\le M}
\frac{\mu_*(j)^2}{j v_j} \Big( \sum_{c\le M/j}
\frac{\Lambda(c)}{c}\frac{\mu_*(c)}{cv_c}\log(Kjc)
\Big)^2\\
&\label{logklog2dfinal}\sim -\hat{h}(0)K\log K \sum_{j\le M}
\frac{\mu_*(j)^2}{j
v_j}\Big(\log(M/j)\log(Kj)+\half\log^2(M/j)\Big)^2.
\end{align}
We combine this with the contribution of the previous calculation.
The sum of (\ref{logklog2dfinal}) and (\ref{logklogrlogdfinal}) is
\begin{align}
&\nonumber\hat{h}(0)K\log K\sum_{j\le M} \frac{\mu_*(j)^2}{j
v_j}\Big(
\frac{3}{2}\log^2(M/j)\log^2(Kj)+\frac{4}{3}\log^3(M/j)\log(Kj)+\frac{1}{4}\log^4(M/j)\Big)\\
&\label{r3}\sim \hat{h}(0)K\log^6K \frac{a^3(6+7a+2a^2)}{12}.
\end{align}

$\bullet$ The contribution of the $\log^3 d_2$ term of
(\ref{logterms}) to (\ref{mainquad}) is
\begin{align}
&\label{logklog3d}\frac{\hat{h}(0)K}{3} \sum_{\al}\sum_{e|\al^2}
\sum_{\substack{r_1,r_2\\(r_1,r_2)=1}}\frac{\sqrt{e}}{\al^2 r_1 r_2}
\Big(\sum_{d|r_1^2r_2^2\al^4/e^2}\frac{\log^3 d_2}{\sqrt{d_1}} \Big)
x_{\al r_1}x_{\al r_2}.
\end{align}
We have
\begin{align}
\label{expansionlog3d}&\sum_{d|r_1^2r_2^2\al^4/e^2}\frac{\log^3
d_2}{\sqrt{d_1}}\\
&\nonumber\indent=\sum_{\substack{i+j+k=3\\i,j,k\ge
0}}a_{ijk}\Big(\sum_{c|r_1}\Lambda_i(c)\tau(r_1/c)\Big)\Big(\sum_{c|r_2}\Lambda_j(c)\tau(r_2/c)\Big)\Big(\sum_{c|\al^2/e}\Lambda_k(c)\tau(\al^2/ec)\Big).
\end{align}
Again the contribution of terms with $k\ge 1$ is $O(K\log^5 K)$.
We have seen that $\sum_{c\le
x}\frac{\Lambda_2(c)\mu_*(c)}{c^2v_c}\ll \log x$. Similarly we
have $\sum_{c\le x}\frac{\Lambda_3(c)\mu_*(c)}{c^2v_c}\ll (\log
x)^2$, by considering the cases in which $c$ is product of one,
two or three distinct primes. Since every term of
(\ref{expansionlog3d}) with $k=0$ must have $i>1$ or $j>1$, we get
that all of (\ref{logklog3d}) is $O(K\log ^{5}K)$.

$\bullet$ The contribution of the $\log r_1r_2 \log^2 d_2$ term of
(\ref{logterms}) to (\ref{mainquad}) is
\begin{align}
&\label{logrlog2d}\frac{-\hat{h}(0)K}{2} \sum_{\al}\sum_{e|\al^2}
\sum_{\substack{r_1,r_2\\(r_1,r_2)=1}}\frac{\sqrt{e}}{\al^2 r_1 r_2}
\Big(\sum_{d|r_1^2r_2^2\al^4/e^2}\frac{\log^2 d_2}{\sqrt{d_1}} \Big)
\log r_1 x_{\al r_1}x_{\al r_2}.
\end{align}
By the same analysis used for the previous terms, this equals
\begin{align}
&\nonumber-\hat{h}(0)K \sum_{\al,\beta}\sum_{e|\al^2}\sum_{b,c}
\frac{e^{\half}\tau(\beta)^2\mu(\beta)\tau(\al^2/e)}{\al^2\beta^{2}}\frac{\Lambda(b)}{b}\frac{\Lambda(c)}{c}\\
&\nonumber\hspace{1in} \times\sum_{r_1,r_2}\frac{x_{\al\beta b
r_1} x_{\al\beta c r_2} \tau(r_1)\tau(r_2)\log(b\beta r_1)}{r_1
r_2}+O(K \log^5 K)\\
&\label{logrlog2d1}\sim -\hat{h}(0)K \sum_j \sum_{b,c} v_j
\frac{\Lambda(b)}{b}\frac{\Lambda(c)}{c} \Big((\log b)
y_{bj}y_{cj}+\sum_d \frac{\Lambda(d)\tau(d)}{d}y_{bdj}y_{cj}\Big).
\end{align}
Using the definition of $y_j$ and partial summation,
(\ref{logrlog2d1}) is
\begin{align}
&\nonumber\sim -\hat{h}(0)K\sum_{j\le M}
\frac{\mu_*(j)^2}{jv_j^2}\Big( \sum_{\substack{b\le M/j\\c\le
M/j}}
\frac{\Lambda(b)^2\Lambda(c)}{bc}\frac{\mu_*(b)\mu_*(c)}{bv_bcv_c}\log(bKj)
\log(cKj)\\
&\nonumber\hspace{1in} + \sum_{\substack{bd\le M/j\\ c\le M/j}}
\frac{\Lambda(b)\Lambda(c)\Lambda(d)\tau(d)}{bcd}\frac{\mu_*(bd)\mu_*(c)}{bdv_{bd}cv_c}\log(bdKj)\log(cKj)\Big)\\
&\nonumber \sim -\hat{h}(0)K\sum_{j\le M}
\frac{\mu_*(j)^2}{jv_j^2}\Big(
\Big(\half\log^2(M/j)\log(Kj)+\frac{1}{3}\log^3(M/j)\Big)\Big(\log(M/j)\log(Kj)\\
&\nonumber \hspace{0.4in} +\half\log^2(M/j)\Big)-\Big(\log^2(M/j)\log(Kj)+\frac{2}{3}\log^3(M/j)\Big)\Big(\log(M/j)\log(Kj)\\
&\nonumber\hspace{3.9in}+\half\log^2(M/j)\Big)\Big)\\
&\label{r4}\sim \hat{h}(0) K \log^6 K \frac{a^4(3+2a)^2}{72}.
\end{align}

\subsection{Remaining terms}
We are left to show that (\ref{otherterms}) is small when
$i_0+i_3+i_4<3$. Since $\log r_1r_2 = \log r_1 + \log r_2$ we have
by symmetry that (\ref{otherterms}) is bounded up to a constant by
\begin{align}
\label{x1} K \log^{i_0}K \sum_{\al}\sum_{e|\al^2}
&\frac{\sqrt{e}\log^{i_1}(\al^2/e)}{\al^2}
\\
&\nonumber \sum_{d|r_1^2r_2^2\al^4/e^2} \frac{\log^{i_2}
d_1\log^{i_3}d_2}{\sqrt{d_1}}\sum_{\substack{r_1,r_2\\(r_1,r_2)=1}}
\frac{x_{\al r_1}x_{\al r_2} \log^{i_4} r_1}{r_1 r_2}.
\end{align}
Using (\ref{use}) we see that this is bounded up to a constant by
\begin{multline}
\label{x2}\sum_{\substack{i_{2,1}+i_{2,2}+i_{2,3}=i_2\\i_{3,1}+i_{3,2}+i_{3,3}=i_3}}
K \log^{i_0}K \sum_{\al}\sum_{e|\al^2}
\frac{\sqrt{e}\log^{i_1}(\al^2/e)}{\al^2}\sum_{d|\al^4/e^2} \frac{\log^{i_{2,3}}d_1 \log^{i_{3,3}}d_2}{\sqrt{d_1}}\\
 \sum_{a_1,b_1,a_2,b_2>0} \frac{
\Lambda_{i_{2,1}}(a_1)\Lambda_{i_{2,2}}(a_2)\Lambda_{i_{3,1}}(b_1)\Lambda_{i_{3,2}}(b_2)}{(a_1a_2)^{\frac{3}{2}}b_1b_2}\\
\sum_{\substack{r_1,r_2\\(a_1b_1r_1,a_2b_2r_2)=1}}
\frac{\tau(r_1)\tau(r_2) \log^{i_4} (a_1b_1r_1)}{r_1 r_2}x_{\al
a_1b_1r_1}x_{\al a_2b_2r_2}.
\end{multline}
Using the Mobius function to detect $(r_1,r_2)=1$ we get
\begin{align}
&\label{negterms}\sum_{\substack{r_1,r_2\\(a_1b_1r_1,a_2b_2r_2)=1}}
\frac{\tau(r_1)\tau(r_2) \log^{i_4} (a_1b_1r_1)}{r_1
r_2}x_{\al a_1b_1r_1}x_{\al a_2b_2r_2}\\
&\nonumber
=\sum_{\beta}\frac{\tau(\beta)^2}{\beta^2}\sum_{\substack{r_1\\(a_2b_2,r_1)=1}}\frac{\log^{i_4}(a_1b_1\beta
r_1)\tau(r_1)}{r_1}x_{\al a_1b_1\beta
r_1}\sum_{\substack{r_2\\(a_1b_1,r_2)=1}}\frac{\tau(r_2)}{r_2}x_{\al
a_2b_2\beta r_2}.
\end{align}
Now we detect $(a_2b_2,r_1)=(a_1b_1,r_2)=1$ using the Mobius
function again. Using the fact that $a_1, b_1, a_2$ and $b_2$ have
at most three prime factors and the bound $y_j\ll \log K$, we have
that (\ref{negterms}) is bounded by a constant multiple of
\begin{align}
&\sum_{\substack{c_1|a_1b_1\\
c_2|a_2b_2}} \sum_{\beta} \frac{\tau(\beta)^2}{\beta^2}
\Big(\log^{i_4} K|y_{\al a_1b_1\beta c_2 }y_{\al a_2b_2\beta c_1}|
+ \sum_{d}\frac{\Lambda_{i_4}(d)}{d}|y_{\al a_1b_1\beta c_2
d}y_{\al a_2b_2\beta c_1}|\Big) \\
&\nonumber \ll \log^{2+i_4} K.
\end{align}
So (\ref{x2}) is
\begin{align}
\label{x3} \ll K\log^{2+i_0+i_3+i_4}K \sum_{\al\le
M}\sum_{e|\al^2}
\frac{\sqrt{e}\log^{i_1}(\al^2/e)}{\al^2}\sum_{d|\al^4/e^2}
\frac{\log^{i_{2}}d_1 }{\sqrt{d_1}}.
\end{align}
The inner sum is
\begin{align}
\frac{1}{\al}\sum_{e|\al^2} \frac{\log^{i_1}e}{\sqrt{e}}
\sum_{d|e^2} \frac{\log^{i_{2}}d_1 }{\sqrt{d_1}} \ll
\frac{1}{\al}\sum_{e|\al^2} e^{-\half+\ep}.
\end{align}
So (\ref{x3}) is $\ll K\log^{3+i_0+i_3+i_4}K$ and we see that the
terms with $i_0+i_3+i_4<3$ contribute only $O(K\log^5 K)$.

\subsection{Completing the evaluation}

We now complete the evaluation of ${\mathcal M}_2$ by adding
together the contributions of each term of (\ref{logterms}). The sum
of (\ref{r1}), (\ref{r2}), (\ref{r3}), and (\ref{r4}) has the
pleasing simplification:
\begin{align}
&\nonumber\hat{h}(0)K\log^6 K\Big(
\frac{a(3+3a+a^2)}{18}+\frac{a^2(2+a)^2}{8}+\frac{a^3(6+7a+2a^2)}{12}+\frac{a^4(3+2a)^2}{72}\Big)\\
&= \hat{h}(0)K\log^6 K \frac{a(1+a)^3(3+3a+a^2)}{18}.
\end{align}

\subsection{Proportion of non-vanishing}

Combining the evaluations of ${\mathcal M}_1$ and ${\mathcal
M}_2$, we get for $a< \half$ that
\begin{align}
\frac{{\mathcal M}_1^2}{{\mathcal M}_2}&\sim
\frac{\Big(\hat{h}(0)K\log^3 K
\frac{a(3+3a+a^2)}{6}\Big)^2}{\hat{h}(0)K\log^6 K
\frac{a(1+a)^3(3+3a+a^2)}{18}}\\
&\nonumber=\Big(1-\frac{1}{(1+a)^3}\Big)\frac{\hat{h}(0)K}{2} \sim
\Big(1-\frac{1}{(1+a)^3}\Big) \sum_{k \equiv 0 \bmod 2}
h\Big(\frac{k-1}{K}\Big)\summ_{f\in H_k} 1.
\end{align}

\vskip0.3in \noindent {\bf Acknowledgements.} I am grateful to
Prof. K. Soundararajan for posing this problem and for many
helpful discussions.

\bibliographystyle{amsplain}

\bibliography{symsquare}

\def\cprime{$'$}
\providecommand{\bysame}{\leavevmode\hbox to3em{\hrulefill}\thinspace}
\providecommand{\MR}{\relax\ifhmode\unskip\space\fi MR }
\providecommand{\MRhref}[2]{%
  \href{http://www.ams.org/mathscinet-getitem?mr=#1}{#2}
}
\providecommand{\href}[2]{#2}
\begin{thebibliography}{10}

\bibitem{blomer}
Valentin Blomer, \emph{On the central value of symmetric square
  {$L$}-functions}, Math. Z. \textbf{260} (2008), no.~4, 755--777.

\bibitem{conrsnai}
J.~B. Conrey and N.~C. Snaith, \emph{Applications of the {$L$}-functions ratios
  conjectures}, Proc. Lond. Math. Soc. (3) \textbf{94} (2007), no.~3, 594--646.

\bibitem{dyson}
Freeman~J. Dyson, \emph{Statistical theory of the energy levels of complex
  systems. {I,II, and III}}, J. Mathematical Phys. \textbf{3} (1962), 140--175.
  \MR{MR0143556 (26 \#1111)}

\bibitem{gelbjac}
Stephen Gelbart and Herv{\'e} Jacquet, \emph{A relation between automorphic
  representations of {${\rm GL}(2)$} and {${\rm GL}(3)$}}, Ann. Sci. \'Ecole
  Norm. Sup. (4) \textbf{11} (1978), no.~4, 471--542.

\bibitem{iwamich}
H.~Iwaniec and P.~Michel, \emph{The second moment of the symmetric square
  {$L$}-functions}, Ann. Acad. Sci. Fenn. Math. \textbf{26} (2001), no.~2,
  465--482.

\bibitem{iwa}
Henryk Iwaniec, \emph{Topics in classical automorphic forms}, Graduate Studies
  in Mathematics, vol.~17, American Mathematical Society, Providence, RI, 1997.
  \MR{MR1474964 (98e:11051)}

\bibitem{ils}
Henryk Iwaniec, Wenzhi Luo, and Peter Sarnak, \emph{Low lying zeros of families
  of {$L$}-functions}, Inst. Hautes \'Etudes Sci. Publ. Math. (2000), no.~91,
  55--131 (2001).

\bibitem{sarniwan2}
Henryk Iwaniec and Peter Sarnak, \emph{The non-vanishing of central values of
  automorphic {$L$}-functions and {L}andau-{S}iegel zeros}, Israel J. Math.
  \textbf{120} (2000), no.~, part A, 155--177.

\bibitem{katzsarn}
Nicholas~M. Katz and Peter Sarnak, \emph{Random matrices, {F}robenius
  eigenvalues, and monodromy}, American Mathematical Society Colloquium
  Publications, vol.~45, American Mathematical Society, Providence, RI, 1999.
  \MR{MR1659828 (2000b:11070)}

\bibitem{kowmic2}
E.~Kowalski and P.~Michel, \emph{The analytic rank of {$J\sb 0(q)$} and zeros
  of automorphic {$L$}-functions}, Duke Math. J. \textbf{100} (1999), no.~3,
  503--542.

\bibitem{kowmic}
\bysame, \emph{A lower bound for the rank of {$J\sb 0(q)$}}, Acta Arith.
  \textbf{94} (2000), no.~4, 303--343. \MR{MR1779946 (2003a:11054)}

\bibitem{lau}
Yuk-Kam Lau, \emph{Non-vanishing of symmetric square {$L$}-functions}, Proc.
  Amer. Math. Soc. \textbf{130} (2002), no.~11, 3133--3139 (electronic).

\bibitem{mont}
H.~L. Montgomery, \emph{The pair correlation of zeros of the zeta function},
  Analytic number theory (Proc. Sympos. Pure Math., Vol. XXIV, St. Louis Univ.,
  St. Louis, Mo., 1972), Amer. Math. Soc., Providence, R.I., 1973,
  pp.~181--193. \MR{MR0337821 (49 \#2590)}

\bibitem{snyder}
A.~E. {\"O}zl{\"u}k and C.~Snyder, \emph{On the distribution of the nontrivial
  zeros of quadratic {$L$}-functions close to the real axis}, Acta Arith.
  \textbf{91} (1999), no.~3, 209--228. \MR{MR1735673 (2001h:11116)}

\bibitem{selb}
Atle Selberg, \emph{Contributions to the theory of the {R}iemann
  zeta-function}, Arch. Math. Naturvid. \textbf{48} (1946), no.~5, 89--155.
  \MR{MR0020594 (8,567e)}

\bibitem{shim}
Goro Shimura, \emph{On the holomorphy of certain {D}irichlet series}, Proc.
  London Math. Soc. (3) \textbf{31} (1975), no.~1, 79--98.

\bibitem{soun}
K.~Soundararajan, \emph{Nonvanishing of quadratic {D}irichlet {$L$}-functions
  at {$s=\frac12$}}, Ann. of Math. (2) \textbf{152} (2000), no.~2, 447--488.

\end{thebibliography}

\end{document}